\documentclass[draft,12pt]{amsart}

\usepackage{tabularx}
\usepackage{hhline}

\newtheorem{theorem}{Theorem}[section]
\newtheorem{proposition}[theorem]{Proposition}
\newtheorem{lemma}[theorem]{Lemma}
\newtheorem{cor}[theorem]{Corollary}

\theoremstyle{definition}
\newtheorem{definition}[theorem]{Definition}

\newtheorem{rem}[theorem]{Remark}
\newtheorem{exam}[theorem]{Example}

\newtheorem{construction}[theorem]{Construction}

\newcommand{\gin}{\ensuremath{\mathrm{gin}}}

\newcommand{\Lex}{\ensuremath{\mathrm{Lex}}}

\newcommand{\reg}{\ensuremath{\mathrm{reg}}}
\newcommand{\Tor}{\ensuremath{\mathrm{Tor}}}

\newcommand{\proj}{ \ensuremath{\mathrm{proj\ dim}}}

\newcommand{\lex}{{\mathrm{lex}}}

\def\cocoa{{\hbox{\rm C\kern-.13em o\kern-.07em C\kern-.13em o\kern-.15em A}}}

\newcommand{\M}{ \ensuremath{\mathcal{M}}}
\newcommand{\N}{ \ensuremath{\mathcal{N}}}
\newcommand{\A}{ \ensuremath{\mathcal{A}}}
\newcommand{\R}{ \ensuremath{\mathcal{R}}}

\newcommand{\dob}[1]{^{\ll {#1}\gg}}
\newcommand{\up}[1]{^{<{#1}>}}
\newcommand{\NN}{\mathbb{N}}
\newcommand{\ZZ}{\mathbb{Z}}

\newcommand{\sqcomp}[1]{_{[{#1}]}}
\newcommand{\sqcompleq}[1]{_{[ \geq {#1}]}}
\newcommand{\SqLex}{\ensuremath{\mathrm{SqLex}}}

\newcommand{\set}{ [\hat{n}]}
\newcommand{\area}{ [\hat{n}] \times \ZZ_{>0}}

\begin{document}

\title
{Hilbert functions of $d$-regular ideals}
\author{Satoshi Murai}
\address{
Department of Pure and Applied Mathematics\\
Graduate School of Information Science and Technology\\
Osaka University\\
Toyonaka, Osaka, 560-0043, Japan
}
\email{s-murai@ist.osaka-u.ac.jp
}
\thanks{The author is supported by JSPS Research Fellowships for Young Scientists}

\maketitle

\begin{abstract}
In the present paper,
we characterize all possible Hilbert functions of graded ideals in a polynomial ring
whose regularity is smaller than or equal to $d$,
where $d$ is a positive integer.
In addition,
we prove the following result which is a generalization of Bigatti, Hulett and Pardue's result:
Let $p \geq 0$ and $d>0$ be integers.
If the base field is a field of characteristic $0$ and
there is a graded ideal $I$ whose projective dimension $\mathrm{proj\ dim}(I)$ is smaller than
or equal to $p$ and whose regularity $\mathrm{reg}(I)$ is smaller than or equal to $d$,
then there exists a monomial ideal $L$ having the maximal graded Betti numbers
among graded ideals $J$ which have the same Hilbert function as $I$ and
which satisfy $\mathrm{proj\ dim}(J) \leq p$ and $\mathrm{reg}(J) \leq d$.
We also prove the same fact for squarefree monomial ideals.
The main methods for proofs are generic initial ideals and combinatorics on strongly stable ideals.
\end{abstract}

\section*{Introduction}
Let $S=K[x_1,\dots,x_n]$ be a standard graded polynomial ring 
over a field $K$.
The (\textit{Castelnuovo--Mumford}) \textit{regularity} of a finitely generated graded $S$-module $M$ is the integer
$ \reg(M)=\max\{k:\beta_{ii+k}(M) \ne 0 \ \mbox{ for some } i \geq 0\},$
where $\beta_{ij}(M)
$ are the graded Betti numbers of $M$. 
In this paper,
we characterize all possible Hilbert functions of graded ideals in $S$
whose regularity is smaller than or equal to $d$ for a given integer $d>0$.

About the characterization of Hilbert functions,
the first important result is the characterization of all possible Hilbert functions of graded ideals in a polynomial ring
which was given by Macaulay \cite{M}.
Another important result is the characterization of all possible Hilbert functions
of squarefree monomial ideals which was given by Kruskal \cite{Kru} and Katona \cite{Ka}
(actually, they characterize face vectors of simplicial complexes).
In the last few decades
this kind of result is particularly of interest in commutative algebra as well as in combinatorics.
In this paper,
we study Hilbert functions of graded ideals $I$ in $S$  with $\reg(I) \leq d$ for a given integer $d >0$.

The basic idea for proofs is generic initial ideals.
A famous result by Bayer and Stillman \cite{BS} implies that,
for any graded ideal $I \subset S$ with $\reg(I) \leq d$,
there exists a strongly stable ideal $J$ such that $I$ and $J$
have the same Hilbert function and $\reg(J) \leq d$.
Thus, to characterize Hilbert functions of graded ideals $I \subset S$ with $\reg(I) \leq d$,
it is enough to consider strongly stable ideals.
Furthermore, the Eliahou--Kervaire formula \cite{EK} says that
a strongly stable ideal $I \subset S$ satisfies $\reg(I) \leq d$ 
if and only if it has no generators of degree $>d$.
Thus characterizing Hilbert functions of graded ideals $I \subset S$
with $\reg(I) \leq d$ is completely a combinatorial problem.
Our results in this paper are based on this idea.
These basic facts on generic initial ideals will be discussed in \S 1.

The vector $(a_1,a_2,\dots,a_n) \in \NN^n$ is called an \textit{$M$-vector}
if there exists an integer $m >0$ and a nonunit graded ideal $I \subset R=k[x_1,\dots,x_m]$
such that
$a_{t}=H(R/I,t-1)$ for $t=1,2,\dots,n$,
where $H(R/I,t)$ is the Hilbert function of $R/I$.
The main result of this paper is the following.

\begin{theorem} \label{main1}
Fix a positive integer $d$.
Let $H: \NN \to \NN$ be a numerical function.
Then there exists a nonzero nonunit graded ideal $I \subset S$
such that $\reg(I) \leq d$ and $H(I,t)=H(t)$ for all $t \in \NN$
if and only if
$H$ satisfies the following conditions.
\begin{itemize}
\item[(i)] There exists a sequence $\ell=(\ell_1,\ell_2,\dots,\ell_n) \in \NN^n$ such that
\begin{itemize}
\item[(a)] $\ell$ is an $M$-vector with $\ell_2 \leq d$;
\item[(b)] $H(t)= \sum_{k=1}^n \ell_k { n-k +t-d \choose n-k}$ for all
$t \geq d$;
\item[(c)] $H(d-1) \leq \ell_n$.
\end{itemize} 
\item[(ii)] $H(0)=0$ and $H(t)\up {n-1} \leq H(t+1)$ for $t < d-1$.
\end{itemize}
\end{theorem}
See Section 3 for the definition of $H(t)\up {n-1}$.
We also characterize all possible Hilbert functions of graded ideals
$I \subset S$ with $\reg(I)=d$ for a fixed integer $d>0$
(Theorem \ref{equalchara}).


A monomial ideal $I \subset S$ is said to be \textit{lexsegment},
if for all monomials $u \in I$ and $v >_\lex u$ with $\deg(v)=\deg(u)$, it follows that $v \in I$,
where $<_\lex$ is the degree lexicographic order induced by $x_1>x_2>\cdots >x_n$.
Macaulay's theorem \cite{M} guarantees that,
for any graded ideal $I \subset S$,
there exists the unique lexsegment ideal $\Lex(I) \subset S$
with the same Hilbert function as $I$.
To prove the main theorem,
we first introduce monomial ideals which play a role similar to lexsegment ideals in the set of
graded ideals $I \subset S$ with $\reg(I) \leq d$ in \S 2.

The complete proof of Theorem \ref{main1} will be given in \S 3.
In this section, we also study the regularity of graded ideals in $S$ with a fixed Hilbert function.
Let $H: \NN \to \NN$ be a numerical function.
Assume that $H$ is the Hilbert function of $S/I$
for some graded ideal $I \subset S$.
Set $a=\min \{\reg(J): H(S/J,t)=H(t) \mbox{ for all }t \in \NN\}$
and $b=\max \{\reg(J): H(S/J,t)=H(t) \mbox{ for all }t \in \NN\}$.
We will show that, for any integer $a \leq r \leq b$,
there exists a monomial ideal $J \subset S$
such that $\reg(J)=r$ and $H(S/J,t)=H(t)$ for all $t \in \NN$.

Let $1 \leq d \leq n$ be an integer.
A simplicial complex $\Gamma$ whose Stanley--Reisner ideal $I_\Gamma$ satisfies $\reg(I_\Gamma) \leq d$
is called a \textit{$(d-1)$-Leray} simplicial complex.
The characterization of face vectors of Leray simplicial complexes was
given by Eckhoff and Kalai.
We refer the reader to \cite{K10} and \cite{K11} for the detail.
Since considering the face vector of a simplicial complex
is equivalent to considering the Hilbert function of its Stanley--Reisner ideal,
their result characterizes Hilbert functions of squarefree monomial ideals $I \subset S$
satisfying $\reg(I) \leq d$.
We will study the result by Eckhoff and Kalai from an algebraic viewpoint in \S 4.
We will show that if $\mathrm{char}(K)=0$ and $H$ is the Hilbert function of some
squarefree monomial ideal $I \subset S$ with $\reg(I) \leq d$,
then there exists a squarefree monomial ideal $L \subset S$
such that $L$ has the maximal graded Betti numbers among squarefree monomial ideals
$ J \subset S$ whose regularity is smaller than or equal to $d$ and whose Hilbert function is equal to $H$.

Bigatti \cite{B}, Hulett \cite{Hu} and Pardue \cite{P}
proved that the graded Betti numbers of a lexsegment ideal
are maximal among graded ideals having the same Hilbert function.
A similar result for squarefree monomial ideals was proved in \cite{AHHshifting} and \cite{MH}.
In particular,
it was proved in \cite{AHHshifting} that if $\mathrm{char}(K)=0$ then
considering the existence of a squarefree monomial ideal having the maximal graded Betti numbers
among squarefree monomial ideals which have the same Hilbert function
is equivalent to considering the existence of a strongly stable ideal having the maximal graded Betti numbers
among graded ideals $I \subset S$ which have the same Hilbert function and
which satisfy $\beta_{ij}(I) =0$ for all $j>n$.
From this viewpoint,
we generalize Bigatti, Hulett and Pardue's result.
We show that there exists a graded ideal having the maximal graded Betti numbers
among graded ideals which have the same Hilbert function and
which accept some restrictions on their graded Betti numbers.
For example,
our results guarantee the existence of a monomial ideal $I \subset S$
having the maximal graded Betti numbers
among graded ideals
which have the same Hilbert function and
which satisfy $\proj(S/I) \leq p$ and $\reg(I) \leq d$ for given positive integers $p$ and $d$,
where $\proj(S/I)=\max\{i:\beta_{ij}(S/I)\ne 0 \mbox{ for some } j\geq 0 \}$ is the projective dimension of $S/I$.
Moreover, we also prove the same fact for squarefree monomial ideals.
These results will be given in \S 5.


\section{Regularity and generic initial ideals}

In this section,
we recall some known results about regularity and generic initial ideals.
Let $S=K[x_1,\dots,x_n]$ be a
standard graded polynomial ring over a field $K$,
$M$ a finitely generated graded $S$-module
and $\beta_{ij}(M)=\dim_K \Tor_i(M,K)_j$ the $i,j$-th graded Betti number of $M$.

If $I$ is a graded ideal in the polynomial ring over an infinite field,
we write $\gin(I)$ for the generic initial ideal of $I$
with respect to the degree reverse lexicographic order
induced by $x_1>x_2>\cdots>x_n$.
We refer the reader to \cite{G} about fundamental facts on generic initial ideals.

A monomial ideal $I \subset S$ is said to be \textit{strongly stable}
if $ux_q \in I$ and $1 \leq p <q$ imply $ux_p \in I$.
For a finitely generated graded $S$-module $M$,
the \textit{Hilbert function} $H(M,t): \mathbb{Z} \to \mathbb{Z}$
of $M$ is the function defined by $H(M,t)=\dim_K M_t$ for $t \in \mathbb{Z}$,
where $M_t$ is the homogeneous component of degree $t$ of $M$.
The following facts are known
(see \cite{G} or \cite{H}).

\begin{lemma} \label{itumono}
Assume that $K$ is infinite.
Let $I \subset S$ be a graded ideal.
\begin{itemize}
\item[(i)]
$I$ and $\gin(I)$ have the same Hilbert function;
\item[(ii)]
If $\mathrm{char}(K)=0$ then $\gin(I)$ is strongly stable;
\item[(iii)]
$\beta_{ij}(I) \leq \beta_{ij}(\gin(I))$ for all $i$ and $j$.
\end{itemize}
\end{lemma}

Let $M$ be a finitely generated graded $S$-module.
The $i,i+j$-th graded Betti number $\beta_{ii+j}(M)$ of $M$
is called an \textit{extremal Betti number of $M$}
if $\beta_{pp+q}(M)=0$ for all $(p,q) \ne (i,j)$ with $p \geq i$
and $q \geq j$.
Bayer, Charalambous and Popescu proved the following nice result.

\begin{lemma}[Bayer--Charalambous--Popescu \cite{BCP}]  \label{bcp}
Assume that $K$ is infinite.
Let $I \subset S$ be a graded ideal.
If $\beta_{ij}(I)$ is an extremal Betti number of $I$
then $\beta_{ij}(\gin(I))$ is an extremal Betti number of $\gin(I)$
and $\beta_{ij}(I)=\beta_{ij}(\gin(I))$.
\end{lemma}
 
\begin{cor} [Bayer--Stillman \cite{BS}]\label{bs}
Assume that $K$ is infinite.
For any graded ideal $I\subset S$,
one has $\reg(I)=\reg(\gin(I))$.
\end{cor}

Another important fact in the theory of generic initial ideals
is that the graded Betti numbers of strongly stable ideals can be
computed by a simple combinatorial method.
For any monomial $u \in  S$,
we write $\max(u)$ for the maximal integer $i$ such that
$x_i$ divides $u$.
Let $I \subset S$ be a monomial ideal.
Then $I$ is said to be \textit{stable}
if $u \in I$ and $k < \max(u)$ imply
$u (x_k / x_{\max(u)}) \in I$.
Clearly, strongly stable ideals are stable.
We write $G(I)$ for the set of minimal monomial generators of $I$. 
The following formula is known as the Eliahou--Kervaire formula.

\begin{lemma}[Eliahou--Kervaire \cite{EK}]\label{ek}
If $I \subset S$ is a stable ideal then
\begin{itemize}
\item[(i)]
$\beta_{ii+k}(I)= \sum_{u \in G(I),\ \deg(u)=k} {\max(u)-1 \choose i}
\ \mbox{ for all }i \mbox{ and }k;$
\item[(ii)] $\reg(I)=\max\{ \deg(u): u \in G(I)\}$.
\end{itemize}
In particular,
the graded Betti numbers of stable ideals are independent
of the characteristic of the base field $K$.
\end{lemma}

Let $d$ be a positive integer and $I\subset S$ a graded ideal.
Then, $I$ is said to be \textit{$d$-regular} if $\reg(I) \leq d$.
Also, we say that $I$ has a \textit{linear resolution} if 
$I$ is generated in degree $d$ and $\reg(I)=d$.
For a positive integer $k$,
we write $I_{\geq k}$ (resp. $I_{\leq k}$) for the ideal generated by all polynomials in $I$ of degree $\geq k$ (resp. $\leq k$).
The following fact is known (see {\cite{EG}}).

\begin{lemma}\label{eg}
A graded ideal $I\subset S$ is $d$-regular
if and only if $I_{\geq d}$ has a linear resolution.
\end{lemma}

If $\mathrm{char}(K) \ne 0$
then $\gin(I)$ is not always strongly stable.
However, the next property easily follows from \cite[Proposition 10]{ERT}.

\begin{lemma} \label{ert}
Assume that $K$ is infinite.
If a graded ideal $I \subset S$ is $d$-regular
then $\gin(I)_{\geq d}$ is a stable ideal generated in degree $d$.
\end{lemma}

\begin{proof}
It follows from \cite[Theorem 1.27]{G} and Corollary \ref{bs} that
$\gin(I)$ is a Borel-fixed monomial ideal with $\reg(\gin(I))=\reg(I)$.
On the other hand, \cite[Proposition 10]{ERT} says that,
if $J \subset S$ is a Borel-fixed monomial ideal then
$J_{\geq \reg(J)}$ is a stable ideal generated in degree $\reg(J)$.
Then the claim follows.
\end{proof}

\begin{cor} \label{dregular}
If $I\subset S$ is a $d$-regular graded ideal then there exists a
strongly stable ideal $J\subset S$ such that $J$ is $d$-regular
and has the same Hilbert function as $I$.
\end{cor}

\begin{proof}
By considering an extension field of $K$,
we may assume that $K$ is infinite.
Then, Lemma \ref{ert} says that
$\gin(I)_{\geq d}$  is a stable ideal generated in degree $d$.
Let $\tilde K$ be a field of characteristic $0$ and $\tilde S =\tilde K [x_1,\dots,x_n]$.
Let $I' \subset \tilde S$ be the monomial ideal generated by $G(\gin(I))$.
Since $I'_{\geq d} \subset \tilde S$ is a stable ideal generated
in degree $d$,
Lemma \ref{ek} says that $I'_{\geq d}$ has a linear  resolution.
Thus $I' \subset \tilde S$ is $d$-regular by Lemma \ref{eg}.
Then the generic initial ideal $\gin(I') \subset \tilde S$ of $I'$
is strongly stable and $d$-regular.
Let $J \subset S$ be the monomial ideal generated by $G(\gin(I'))$.
Clearly, $J$ is a strongly stable ideal.
Since the graded Betti numbers of strongly stable ideals are independent of the characteristic of a base field,
the ideal $J \subset S$ is $d$-regular and has the same Hilbert function as $I$
by Lemma \ref{itumono} (i).
%
\end{proof}

\section{Combinatorics on strongly stable ideals}

In this section,
we will study strongly stable ideals to prove the main theorem.
Macaulay's characterization of Hilbert functions of graded ideals
says that, for any graded ideal $I \subset S$,
there exists the unique lexsegment ideal with the same Hilbert function as $I$.
The aim of this section is to introduce strongly stable ideals which
play a role similar to lexsegment ideals in the set of $d$-regular graded ideals.

Let $\M_d$ be the set of all monomials in $S=K[x_1,\dots,x_n]$ of degree $d$.
A set of monomials $V \subset S$ is said to be
\textit{strongly stable} if $ux_q \in V$ and $1\leq p<q$ imply $ux_p \in V$.
Also,
a set of monomials $L \subset S$ is said to be \textit{lexsegment} 
if, for all monomials $u \in L$ and $v >_\lex u$ with $\deg(v)=\deg(u)$, it follows that $v \in L$.
Note that if $L \subset S$ is a lexsegment set of monomials
then $L$ is strongly stable.

For a subset $V \subset \M_d$,
we write
$$M_{\leq k}(V)=\{u \in V: \max(u) \leq k\} \ \ \
\mbox{ for } k=1,2,\dots,n.$$
For a finite set $A$, let $|A|$ denote the cardinality of $A$.
We recall the following fact which is due to Bayer.
(See {\cite[\S 1]{B}}.)

\begin{lemma}\label{bigatti}
Let $L \subset \M_d$ be a lexsegment,
and let $V \subset \M_d$ be a strongly stable set of monomials with $|L|=|V|$.
Then
\begin{itemize}
\item[(i)] $|M_{\leq k}(V)| \geq |M_{\leq k}(L)|  \mbox{ for } k=1,2,\dots,n$;
\item[(ii)]
$M_{\leq k}(L)$ is also lexsegment in $K[x_1,\dots,x_k]$ for $k=1,2,\dots,n$.
\end{itemize}
\end{lemma}

\begin{definition}
Let $V \subset \M_d$.
Then $V$ can be written in the form.
\begin{eqnarray}
&&V = \bigcup_{k=1}^n\{u \in V: \max(u)=k\} \label{SRM1}
 = \bigcup_{k=1}^n x_k \{(u/ x_k): u \in V \mbox{ and }\max(u)=k\}.
\end{eqnarray}
Set
$$D_k(V)=\{(u/ x_k): u \in V \mbox{ and }\max(u)=k\} \ \ \mbox{ for }k=1,2,\dots,n.$$
Note that each $D_k(V)$ is a set of monomials in $K[x_1,\dots,x_k]$ of degree $d-1$.
Also, if $V$ is strongly stable then each $D_k(V)$ is strongly stable.

For a strongly stable set of monomials $V \subset \M_d$,
we define the sequence $\ell(V)=(\ell_1(V),\ell_2(V),\dots,\ell_n(V))$
by $\ell_k(V)=|D_k(V)|$ for $k=1,2,\dots,n$.
In other words, $\ell_k(V)$ is the number of monomials $u \in V$ with $\max(u) =k$.
If $I$ is the monomial ideal generated by a strongly stable set of monomials $V \subset \M_d$,
then we write $\ell(I)=\ell(V)$.
The sequence $\ell(V)$ (resp.\ $\ell(I)$)
is called the \textit{$\ell$-sequence of} $V$ (resp.\ $I$).
\end{definition}

Let $I$ be a strongly stable ideal generated in degree $d$.
Then the Eliahou--Kervaire formula (Lemma \ref{ek}) says that $I$ has a $d$-linear resolution.
It is not hard to see that the $\ell$-sequence of $I$
determines the Hilbert function of $I$.
Indeed, by using the Eliahou--Kervaire formula, we have
\begin{eqnarray}
\beta_{ii+d}(I)= \sum_{k=1}^n \ell_k(I) {k-1 \choose i} 
\ \ \mbox{ for all }i. \label{LRM2}
\end{eqnarray}
The above formula implies the following fact.

\begin{lemma} \label{equive}
Let $I \subset S$ and $I' \subset S$ be strongly stable ideals generated in degree $d$.
Then the following conditions are equivalent.
\begin{itemize}
\item[(i)]  $I$ and $I'$ have the same Hilbert function;
\item[(ii)]  $I$ and $I'$ have the same graded Betti numbers;
\item[(iii)]  $\ell(I)=\ell(I')$.
\end{itemize}
\end{lemma}

\begin{proof}
Since $I$ and $I'$ have a $d$-linear resolution,
(i) $\Leftrightarrow$ (ii) is obvious.
On the other hand, (ii) $\Leftrightarrow$ (iii)
follows from (\ref{LRM2}).
\end{proof}

\begin{exam} \label{ex1}
Let $V=\{x_1^3,x_1^2x_2,x_1x_2^2,x_2^3,x_1^2x_3,x_1x_2x_3,x_2^2x_3,x_1^2x_4\}.$
Then
\begin{eqnarray*}
D_1(V)&=&\{x_1^2\}, \\
D_2(V)&=&\{x_1^2,x_1x_2,x_2^2\}, \\
D_3(V)&=&\{x_1^2,x_1x_2,x_2^2\}, \\
D_4(V)&=&\{x_1^2\}
\end{eqnarray*}
and $\ell(V)=(1,3,3,1)$.
\end{exam}

\begin{definition}
We say that a set of monomials $V \subset \M_d$
is \textit{$d$-linear lexsegment}
if $V$ is strongly stable and each $D_k(V) \subset \M_{d-1}$
is lexsegment in $K[x_1,\dots,x_k]$ for $k=1,2,\dots,n$. 
\end{definition}

\begin{exam}
It is easily verified that if $V \subset \M_d$ is lexsegment
then $V$ is $d$-linear lexsegment.
Also, any $d$-linear lexsegment subset of $\M_d$ is strongly stable.
We give examples which show that these three classes are different.

Set $L=\{x_1^3,x_1^2x_2,x_1x_2^2,x_2^3,x_1^2x_3\}$.
Then $L$ is strongly stable and $D_1(L)=\{x_1^2\}$, $D_2(L)=\{x_1^2,x_1x_2,x_2^2\}$
and $D_3(L)=\{x_1^2\}$.
Each $D_k(L)$ is a lexsegment set of monomials in $K[x_1,\dots,x_k]$.
Thus this set $L$ is $3$-linear lexsegment, however,
is not lexsegment since $x_1x_2x_3>_\lex x_2^3$.
On the other hand,
the set of monomials $V$ given in Example \ref{ex1}
is strongly stable but is not $3$-linear lexsegment since
$D_3(V)$ is not lexsegment in $K[x_1,x_2,x_3]$.
\end{exam}

The idea of $d$-linear lexsegment subsets is useful
to study the $\ell$-sequence of strongly stable
subsets $V\subset \M_d$. 
Indeed, we have

\begin{proposition} \label{llex}
If $V$ is a strongly stable subset of $\M_d$
then there exists the unique $d$-linear lexsegment subset $L \subset \M_d$
such that $\ell(L)=\ell(V)$.
\end{proposition}

To prove the above proposition,
we first show the next lemma.

\begin{lemma} \label{hodai}
Let $V \subset \M_d$.
Then $V$ is strongly stable if and only if $V$ satisfies
\begin{itemize}
\item[(i)] each $D_k(V)$ is strongly stable for $k=1,2,\dots,n$;
\item[(ii)] 
$M_{\leq k-1}(D_k(V)) \subset D_{k-1}(V)$
for $k=2,3,\dots,n$. 
\end{itemize}
\end{lemma}

\begin{proof}
(``Only if")
Assume that $V$ is strongly stable.
It is clear that $V$ satisfies condition (i).
We will consider condition (ii).

Let $u \in M_{\leq k-1}(D_k(V))$.
Then $\max(u) \leq k-1$ and $ux_k \in V$.
Since $V$ is strongly stable, we have $u x_{k-1} \in V$ and $\max(ux_{k-1})= k-1$.
Hence we have $u \in D_{k-1}(V)$.
Thus we have $M_{\leq k-1}(D_k(V)) \subset D_{k-1}(V)$ as desired.

(``If")
Assume that $V$ satisfies condition (i) and (ii).
We will show that $V$ is strongly stable.
Let $u \in V$ and $k=\max(u)$.
Write $u=v x_k$.
Let $1 \leq s <t$ be integers such that $x_t$ divides $u$.
What we must prove is $u(x_s/x_t) \in V$.

If $x_t$ divides $v$ then $u(x_s /x_t)= \{v (x_s/x_t)\}x_k \in x_k D_k(V)\subset V$
by condition (i).

Assume that $x_t$ does not divide $v$.
Then $t=k$ and $\max(v)<k$.
Also, condition (ii) says that, for any integer $1 \leq p<k$,
one has
\begin{eqnarray}
\left.
\begin{array} {lll}
D_k(V) \cap K[x_1,\dots,x_p] & = &
(M_{\leq k-1}(D_k(V))) \cap K[x_1,\dots,x_p]  \medskip\\
& \subset & D_{k-1}(V) \cap K[x_1,\dots,x_p].
\end{array}
\right.
\label{LRM3}
\end{eqnarray}
Set $l=\max(v)<k$.
Then $v \in D_k(V) \cap K[x_1,\dots,x_l]$.
Also (\ref{LRM3}) says that
\begin{eqnarray*}
\ D_k(V) \cap K[x_1,\dots,x_l] \subset D_{k-1}(V) \cap K[x_1,\dots,x_l]
 \subset \cdots \subset D_l(V) \cap K[x_1,\dots,x_l].
\end{eqnarray*}
Thus we have $v \in D_p(V)$ for all $l \leq p \leq k$.
Recall that $u(x_s/x_t)= v x_s$ since $t=k$.
If $s \geq l$ then the above fact says $u(x_s/x_t)= v x_s \in x_s D_s(V) \subset V$.
On the other hand,
if $s<l$ then, since $v \in D_l(V)$ and $D_l(V)$ is strongly stable by condition (i),
we have
$$ u(x_s/x_t) = v x_s =\{ v (x_s/x_l) \} x_l \in x_l D_l(V) \subset V.$$
In both cases, we have $u(x_s/x_t) \in V$.
Hence $V$ is strongly stable.
\end{proof}

\begin{proof}[Proof of Proposition \ref{llex}]
Since $d$-linear lexsegment subsets are uniquely determined from their $\ell$-sequence,
it is enough to prove the existence of a subset which satisfies the required conditions.
For $k=1,2,\dots,n$,
let $B_k\subset \M_{d-1} \cap K[x_1,\dots,x_k]$ be the set of monomials
satisfying that $B_k$ is lexsegment in $K[x_1,\dots,x_k]$
and $|B_k|=|D_k(V)|$.
Set $L= \bigcup_{k=1}^n x_k B_k \subset \M_d$.
We will show that this set $L$ is a $d$-linear lexsegment set of monomials with $\ell (L)=\ell(V)$.
Since $D_k(L)=B_k$ is lexsegment in $K[x_1,\dots,x_k]$
and $|B_k|=|D_k(V)|=\ell_k(V)$ for all $k$,
what we must prove is that $L$ is strongly stable.
We will show that $L$ satisfies condition (i) and (ii)
of Lemma \ref{hodai}.

Any lexsegment set of monomials is strongly stable.
Thus $D_k(L)=B_k$ is strongly stable for all $k$.
Hence $L$ satisfies condition (i) of Lemma \ref{hodai}.
On the other hand,
Lemma \ref{bigatti} says that
\begin{eqnarray*}
&&|M_{\leq k-1}(D_k(V))| \geq |M_{\leq k-1}(B_k)| \ \ \mbox{ for }k=2,3,\dots,n.
\end{eqnarray*}
Also, since $V$ is strongly stable,
Lemma \ref{hodai} says that $D_{k-1}(V) \supset M_{\leq k-1}(D_k(V))$.
Then we have
\begin{eqnarray}
&&|B_{k-1}|=|D_{k-1}(V)| \geq |M_{\leq k-1}(D_k(V))| \geq |M_{\leq k-1}(B_k)| 
\ \ \mbox{ for }k=2,3,\dots,n.  \label{past4}
\end{eqnarray}
Since both $B_{k-1} \subset \M_{d-1}$ and $M_{\leq k-1}(B_k) \subset \M_{d-1}$
are lexsegment in $K[x_1,\dots,x_{k-1}]$ by Lemma \ref{bigatti} (ii),
we have $B_{k-1} \supset M_{\leq k-1}(B_k)$ by (\ref{past4}).
Since $D_k(L)=B_k$ for all $k$,
it follows that $L$ satisfies condition (ii) of Lemma \ref{hodai}.
Thus $L$ is strongly stable as required.
\end{proof}


\begin{definition}
Let $d$ be a positive integer.
We say that a strongly stable ideal $I$ is \textit{$d$-linear lexsegment}
if $I$ is generated by a $d$-linear lexsegment set of monomials $V \subset \M_d$.

For a positive integer $k$ and a monomial ideal $I \subset S$,
we say that $I_k$ is lexsegment 
if the set of all monomials in $I_k$ is lexsegment.
A monomial ideal $I$ is said to be \textit{$d$-lexsegment}
if $I_{\geq d}$ is $d$-linear lexsegment and
$I_k$ is lexsegment for all $k<d$.

In the next section,
it will be shown that, for any $d$-regular graded ideal $I \subset S$,
there exists the unique $d$-lexsegment ideal with the same Hilbert function as $I$.
In this section, we note some easy properties of $d$-lexsegment ideals.
\end{definition}

\begin{proposition} \label{dbasic}
Let $d$ be a positive integer.
\begin{itemize}
\item[(i)] Any $d$-lexsegment ideal $I \subset S$ is strongly stable and $d$-regular;
\item[(ii)] If $I \subset S$ and $J \subset S$ are $d$-lexsegment and
have the same Hilbert function then $I=J$;
\item[(iii)] If $I \subset S$ is a $d$-lexsegment ideal with $\reg(I)<d$
then $I$ is lexsegment.
\end{itemize}
\end{proposition}

\begin{proof}
(i) It is clear that $I$ is strongly stable since lexsegment ideals and $d$-linear lexsegment ideals are strongly stable.
Also, since $I_{\geq d}$ is generated in degree $d$,
$I$ has no generators of degree $>d$.
Then Lemma \ref{ek} says that $I$ is $d$-regular.

(ii) Since $I_k$ and $J_k$ are lexsegment for $k<d$
and $I$ and $J$ have the same Hilbert function, it follows that $I_k=J_k$ for $k<d$.
On the other hand, $I_{\geq d}$ and $J_{ \geq d}$ are $d$-linear lexsegment and have the same Hilbert function.
Then Lemma \ref{equive} says that $\ell(I_{\geq d})=\ell (J_{ \geq d})$.
Since $d$-linear lexsegment ideals are uniquely determined from their $\ell$-sequence, we have $I_{\geq d}=J_{ \geq d}$.

(iii) By the definition of $d$-lexsegment ideals,
it follows that $I_{\leq d-1}$ is lexsegment.
Since $I$ is strongly stable, $\reg(I)<d$ implies
that $I$ has no generators of degree $\geq d$ by Lemma \ref{ek}.
Then $I = I_{\leq d-1}$ is lexsegment.
\end{proof}

\section{The Hilbert function of $d$-regular graded ideals}

In this section,
we will prove Theorem \ref{main1}.
To prove the theorem,
we first characterize
all possible $\ell$-sequences of strongly stable ideals generated in degree $d$
by using combinatorics on Stanley--Reisner ideals.

Before considering Stanley--Reisner ideals,
we recall Macaulay's characterization of Hilbert functions of graded ideals.
Let $d$ be a positive integer.
Any integer $a>0$ can be written uniquely in the form
$$a={ a(d) \choose d} + { a(d-1) \choose d-1}+ 
\cdots +{a(j) \choose j},$$
where $a(d)>a(d-1)> \cdots >a(j) \geq j \geq 1$
(see \cite[Lemma 4.2.6]{CM}).
The above representation is called the \textit{$d$-th
Macaulay representation of $a$.}
If $a={ a(d) \choose d} + { a(d-1) \choose d-1}+ 
\cdots + {a(j) \choose j}$ is the $d$-th Macaulay representation of an integer $a>0$,
then we write
$$a \dob d =
{ a(d)+1 \choose d+1} + { a(d-1)+1 \choose d}+ 
\cdots + {a(j)+1 \choose j+1}$$
and
$$a \up d =
{ a(d)+1 \choose d} + { a(d-1)+1 \choose d-1}+ 
\cdots + {a(j)+1 \choose j}.$$
Also, we set $0 \dob d =0$ and $0 \up d =0$.

\begin{lemma}[Macaulay] \label{mac}
Let $H: \mathbb{N} \to \mathbb{N}$ be a numerical function.
Then
\begin{itemize}
\item[(i)] there exists a nonunit graded ideal $I \subset S$
such that $H(S/I,t)=H(t)$ for all $t \in \NN$
if and only if $H(0)=1$, $H(1) \leq n$ and
$H(t) \dob t \geq H(t+1)$  for all $t \geq 1$.
\item[(ii)] there exists a nonunit graded ideal $I \subset S$
such that $H(I,t)=H(t)$ for all $t \in \NN$
if and only if $H(0)=0$ and
$H(t) \up {n-1} \leq H(t+1) \leq H(S,t+1)$  for all $t \geq 0$.
\item[(iii)] for any graded ideal $I \subset S$,
there exists the unique lexsegment ideal, denoted $\Lex(I) \subset S$,
such that $I$ and $\Lex(I)$ have the same Hilbert function.
\end{itemize}
\end{lemma}

We refer the reader to \cite[Theorem 4.2.10]{CM} for the proof of the
above lemma.
(Statement (ii) was not written in \cite{CM}.
However, this statement is easily verified by using \cite[Proposition 4.2]{B}
together with statement (iii).)

Next, we recall some techniques which were developed in the theory of face vectors of simplicial complexes.
Let $\Gamma$ be a simplicial complex on $[n]=\{1,2,\dots,n\}$,
that is,
$\Gamma$ is a collection of subsets of $[n]$
satisfying that $F \in \Gamma$ and $G \subset F$ implies $G \in \Gamma$
(we do not assume that $\{i\} \in \Gamma$ for all $i \in [n]$).
The \textit{dimension of $\Gamma$} is 
$\dim \Gamma = \max\{|F|: F \in \Gamma\} -1$.
For an integer $k \geq 0$,
define $\Gamma_k=\{ F \in \Gamma: |F|=k+1\}$ and $f_k(\Gamma)=|\Gamma_k|$.
Set $d= \dim \Gamma +1$.
The vector $f(\Gamma)=(f_0(\Gamma),f_1(\Gamma),\dots,f_{d-1}(\Gamma))$
is called the \textit{$f$-vector of $\Gamma$}.
The \textit{$h$-vector} $h(\Gamma)=(h_0(\Gamma),h_1(\Gamma),\dots,h_d(\Gamma))$ of $\Gamma$
is defined by the relation
$$h_k(\Gamma)= \sum_{i=0}^k (-1)^{k-i} {d-i \choose k-i} f_{i-1}(\Gamma),$$
where we let $f_{-1}(\Gamma)=1$.

The \textit{Stanley--Reisner} ideal $I_\Gamma \subset S$ of a simplicial complex $\Gamma$ on $[n]$
is the monomial ideal generated by all squarefree monomials
$x_{i_1}x_{i_2}\cdots x_{i_k} \in S$ with $\{i_1,i_2,\dots,i_k\} \not \in \Gamma$.
We say that a simplicial complex $\Gamma$ on $[n]$ is Cohen--Macaulay
if $S/I_\Gamma$ is Cohen--Macaulay.

A simplicial complex $\Gamma$ is said to be
\textit{shifted} if $F \in \Gamma$ and $j \in F$
implies $(F\setminus \{j\}) \cup \{i\} \in \Gamma$ for any $j<i \leq n$ with $i \not \in F$.
The following fact appeared in the theory of exterior algebraic shifting.
See \cite[\S 7 and \S 8]{H}.

\begin{lemma}[Kalai] \label{shifting}
If $\Gamma$ is a (Cohen--Macaulay) simplicial complex on $[n]$
then there exists a (Cohen--Macaulay) shifted simplicial complex
$\Delta$ on $[n]$ such that $f(\Gamma)=f(\Delta)$ and $\reg(I_\Gamma)=\reg(I_\Delta)$.
\end{lemma}

We recall Stanley's characterization of $h$-vectors of Cohen--Macaulay
simplicial complexes.
A finite sequence of integers $h=(h_0,h_1,\dots,h_p) \in \mathbb{N}^{p+1}$
is called an \textit{$M$-vector}
if $h_0=1$ and $h_t \dob t \geq h_{t+1}$ for all $t \geq 1$.
This definition is the same as the definition given in the introduction by Lemma \ref{mac}.

\begin{lemma}[{Stanley \cite[II Theorem 3.3]{S}}] \label{stanley}
Fix an integer $n \geq d \geq 1$.
Let $h=(h_0,h_1,\dots,h_d) \in \mathbb{N}^{d+1}$.
The following conditions are equivalent.
\begin{itemize}
\item[(a)]
There exists a $(d-1)$-dimensional Cohen--Macaulay simplicial complex $\Gamma$
on $[n]$ such that $h(\Gamma)=h$;
\item[(b)]
There exists a $(d-1)$-dimensional Cohen--Macaulay shifted simplicial complex $\Gamma$
on $[n]$ such that $h(\Gamma)=h$;
\item[(c)]
$h$ is an $M$-vector with $h_1 \leq n-d$.
\end{itemize}
\end{lemma}

Note that the equivalence of (a) and (b) follows from
Lemma \ref{shifting}.
\medskip

A set of squarefree monomials $V \subset S$ is said to be \textit{squarefree strongly stable} if
for all squarefree monomials $u=x_{i_1}\cdots x_{i_k} \in V$, it follows that
$u(x_p/x_q) \in V$ for all integers $1 \leq p<q$ with $p \not \in \{i_1,\dots,i_k\}$ and $q \in  \{i_1,\dots,i_k\}$.
We say that a squarefree monomial ideal $I$ is \textit{squarefree strongly stable}
if the set of squarefree monomials in $I$ is squarefree strongly stable. 
Note that a simplicial complex $\Gamma$ is shifted if and only if
$I_\Gamma$ is squarefree strongly stable.
The graded Betti numbers of squarefree strongly stable ideals can be computed like strongly stable ideals.

\begin{lemma}[Aramova--Herzog--Hibi \cite{SQlex}] \label{sqek}
Let $I \subset S$ be a squarefree strongly stable ideal.
Then
\begin{itemize}
\item[(i)] $\beta_{ii+k}(I)=\sum_{u \in G(I),\ \deg(u)=k}
{\max(u)-k \choose i}$ for all $i$ and $k$;
\item[(ii)] $\reg(I)=\max\{ \deg(u): u \in G(I)\}$.
\end{itemize}
\end{lemma}

Let
$\N_d$ be the set of all squarefree monomials in
$K[x_1,x_2,\dots,x_{n+d-1}]$ of degree $d$.
The \textit{squarefree operation} $\Phi: \bigcup_{d \geq 0} \M_d
\to \bigcup_{d \geq 0} \N_d$ is the map defined by
$$\Phi(x_{i_1} x_{i_2} \cdots x_{i_d})=
x_{i_1} x_{i_2+1} x_{i_3 +2} \cdots x_{i_d + d-1},$$
where $i_1 \leq i_2 \leq \cdots \leq i_d$.
Since the squarefree operation is bijective,
we can define the inverse map 
$\Phi^{-1}:\bigcup_{d \geq 0} \N_d \to \bigcup_{d \geq 0} \M_d$.

Let $I\subset S$ be a strongly stable ideal generated in degree $d$.
We write $\Phi(I)$ for the ideal in $K[x_1,x_2,\dots,x_{n+d-1}]$
generated by all squarefree monomials $\Phi(u)$ with $u \in G(I)$.
Conversely, if $J \subset K[x_1,\dots,x_{n+d-1}]$ is a squarefree strongly stable ideal generated in degree $d$
then we define $\Phi^{-1}(J) \subset S$ in the same way.
The following fact is known.

\begin{lemma}[Aramova--Herzog--Hibi {\cite{AHHshifting}}] \label{operation}
Let $d$ be a positive integer.
If $I \subset S$ is a strongly stable ideal generated in degree $d$,
then $\Phi(I) \subset K[x_1,\dots,x_{n+d-1}]$ is squarefree strongly stable.
Conversely,
if $J\subset K[x_1,\dots,x_{n+d-1}]$ is a squarefree strongly stable ideal generated in degree $d$,
then $\Phi^{-1}(J) \subset S$ is strongly stable.
\end{lemma}

Let $I \subset K[x_1,x_2,\dots,x_{n+d-1}]$ be a squarefree strongly stable ideal generated in degree $d$. 
For $k=1,2,\dots,n$,
let $\ell_k^*(I)$ denote the number of squarefree monomials
$u \in G(I)$ with $\max(u)=k+d-1$.
The sequence $\ell^*(I)=(\ell_1^*(I),\ell_2^*(I),\dots,\ell_n^*(I))$
will be called the $\ell^*$-sequence of $I$.
It is easily verified that $\ell^*(I)=\ell(\Phi^{-1}(I))$.

The \textit{Alexander dual} of a simplicial complex $\Gamma$ on $[n]$
is the simplicial complex on $[n]$ defined by
$$\Gamma^* =\{ F \subset [n]: [n]\setminus F \not \in \Gamma\}.$$

\begin{lemma}[{Eagon--Reiner \cite[Theorems 3 and 4]{ER}}] \label{er}
Let $d$ be a positive integer
and $\Gamma$ an $(n-2)$-dimensional simplicial complex on $[n+d-1]$.
Then $\Gamma$ is Cohen--Macaulay if and only if
$I_{\Gamma^*}$ has a $d$-linear resolution.
Moreover, if $\Gamma$ is Cohen--Macaulay then
$$\beta_{ii+d}(I_{\Gamma^*})= \sum_{k=0}^{n-1}  h_k(\Gamma) { k \choose i} \ \ \mbox{ for all i}.$$
\end{lemma}

\begin{cor} \label{landh}
Let $d$ be a positive integer and
$\Gamma$ an $(n-2)$-dimensional simplicial complex on $[n+d-1]$.
Then $\Gamma$ is Cohen--Macaulay and shifted if and only if
$I_{\Gamma^*}$ is a squarefree strongly stable ideal generated in degree $d$.
Moreover, if $\Gamma$ is Cohen--Macaulay and shifted then $\ell^*(I_{\Gamma^*})=h(\Gamma)$.
\end{cor}

\begin{proof}
It is easy to see that $\Gamma$ is shifted if and only if $\Gamma^*$ is shifted.
Then the first statement follows from Lemmas \ref{sqek} and \ref{er}.
Also, by Lemma \ref{sqek}, if $I_{\Gamma^*}$ is a squarefree strongly stable ideal generated in degree $d$
then
$$\beta_{ii+d}(I_{\Gamma^*})= \sum_{k=1}^n  \ell^*_k(I_{\Gamma^*}) { k-1 \choose i} \ \ \mbox{ for all }i.$$
The above formula together with Lemma \ref{er}
implies $\ell^*(I_{\Gamma^*})=h(\Gamma)$.
\end{proof}

Now, we are in the position to characterize all possible $\ell$-sequences
of strongly stable ideals generated in degree $d$.

\begin{proposition} \label{lsequence}
Fix a positive integer $d$.
Let $\ell=(\ell_1,\ell_2,\dots,\ell_n) \in \mathbb{N}^n$.
The following conditions are equivalent.
\begin{itemize}
\item[(a)] There exists a nonzero strongly stable ideal $I \subset K[x_1,\dots,x_n]$
generated in degree $d$ such that $\ell(I)=\ell$;
\item[(b)] There exists a nonzero squarefree strongly stable ideal $I \subset
K[x_1,\dots,x_{n+d-1}]$ generated in degree $d$ such that $\ell^*(I)=\ell$;
\item[(c)] $\ell$ is an $M$-vector with $\ell_2 \leq d$.
\end{itemize}
\end{proposition}

\begin{proof}
The equivalence of (a) and (b) follows from Lemma \ref{operation}
together with the fact that $\ell(I)=\ell^*(\Phi(I))$ for any strongly stable ideal $I \subset S$ generated in degree $d$.
On the other hand,
Corollary \ref{landh} says that (b) is equivalent to the condition that
there exists an $(n-2)$-dimensional Cohen--Macaulay shifted simplicial complex $\Gamma$ on $[n+d-1]$
such that $h(\Gamma)=\ell$.
Then the equivalence of (b) and (c) follows from Lemma \ref{stanley}.
\end{proof}

By Lemma \ref{equive},
if $I\subset S$ is a strongly stable ideal generated in degree $d$,
then the $\ell$-sequence of $I$ must determine the Hilbert function of $I$.
We write how to determine the Hilbert function of $I$ from $\ell(I)$.

\begin{lemma} \label{hilb}
Let $I\subset S$ be a strongly stable ideal generated in degree $d$.
Then
$$H(I,t)=\sum_{k=1}^n \ell_k (I) {n-k+t-d \choose n-k}.$$
\end{lemma}

\begin{proof}
Let $G(I) =\{u_1,\dots,u_t\}$.
It follows from \cite[Lemma 1.1]{EK} that
any monomial $u \in I$ can be written uniquely in the form
$u=v w$ with $v \in G(I)$ and $\max(v) \leq \min (w)$,
where $\min (w)$ is the minimal integer $p$ such that $x_p$ divides $w$.
This fact says that $I$ can be written in the form
$$I= \bigoplus_{k=1}^t u_k K[x_{\max(u_k)},\dots,x_n].$$
Since $H(K[x_1,\dots,x_n],t)= { n-1 +t \choose n-1}$,
the above equation implies the claim.
\end{proof}

\begin{cor}
Let $d$ be a positive integer and $H: \NN \to \NN$ a numerical function.
Then there exists a graded ideal $I\subset S$ such that
$I$ has a $d$-linear resolution and
$H(I,t)=H(t)$ for all $t \in \NN$
if and only if
there exists a sequence $\ell=(\ell_1,\ell_2,\dots,\ell_n) \in \NN^n$ such that
$\ell$ is an $M$-vector with $\ell_2 \leq d$ and
$H(t)= \sum_{k=1}^n \ell_k { n-k +t-d \choose n-k}$ for all $t \in \NN$.
\end{cor}

\begin{proof}
By Corollary \ref{dregular}, for any graded ideal $I\subset S$ which has a $d$-linear resolution,
there exists a strongly stable ideal $J \subset S$ such that
$J$ has a $d$-linear resolution and has the same Hilbert function as $I$.
Then the claim follows from Proposition \ref{lsequence} and Lemma \ref{hilb}.
\end{proof}

Now, we give a proof of the necessity of Theorem \ref{main1}.

\begin{theorem} \label{ness}
Fix a positive integer $d$.
Let $I \subset S$ be a nonzero nonunit $d$-regular graded ideal,
and let $H: \NN \to \NN$ be the Hilbert function of $I$.
Then
\begin{itemize}
\item[(i)] there exists a sequence $\ell=(\ell_1,\ell_2,\dots,\ell_n) \in \NN^n$ such that
\begin{itemize}
\item[(a)] $\ell$ is an $M$-vector with $\ell_2 \leq d$;
\item[(b)] $H(t)= \sum_{k=1}^n \ell_k { n-k +t-d \choose n-k}$ for all
$t \geq d$;
\item[(c)] $H(d-1) \leq \ell_n$.
\end{itemize} 
\item[(ii)] $H(0)=0$ and $H(t)\up {n-1} \leq H(t+1)$ for $t < d-1$.
\end{itemize}
\end{theorem}

\begin{proof}
It is clear that $H$ satisfies condition (ii) by Lemma \ref{mac}.
We will consider condition (i).

We may assume that $I$ is strongly stable by Corollary \ref{dregular}.
Since $I$ is $d$-regular,
$I_{\geq d}$ is a strongly stable ideal generated in degree $d$
by Lemma \ref{ek}.
Set  $\ell=\ell(I_{\geq d})$.
Then $\ell$ satisfies condition (i)--(a) and (i)--(b)
by Proposition \ref{lsequence} and Lemma \ref{hilb}.
Set $A=\{ u \in \M_{d-1}: u \in I\}$
and $V=G(I_{\geq d})=\{ u \in \M_d: u \in I\}$.
Since $x_n A \subset V= \bigcup_{k=1}^n x_k D_k(V)$,
we have $A \subset D_n(V)$.
Hence we have $H(I,d-1)=|A| \leq |D_n(V)| = \ell_n(I_{\geq d})$.
Thus $\ell$ also satisfies condition (i)--(c).
\end{proof}


Next, we will give a proof of the sufficiency of Theorem \ref{main1}.

\begin{theorem} \label{suff}
Let $H: \NN \to \NN$ be a numerical function satisfying condition 
$\mathrm{(i)}$ and $\mathrm{(ii)}$ of Theorem \ref{ness}.
Then there exists a $d$-lexsegment ideal $J \subset S$ 
such that $H(J,t)=H(t)$ for all $t \in \NN.$
\end{theorem}

\begin{proof}
Condition (i)--(a) and (i)--(c) say $H(d-1) \leq { d+n-2 \choose n-1} =H(S,d-1)$.
Then, inductively,
condition (ii) says $H(k) \leq { k +n-1 \choose n-1} = H(S,k) $ for $k \leq d-1$.
%
%
Thus Lemma \ref{mac} says that
there exists the lexsegment ideal $I \subset S$
such that $H(I,t)=H(t)$ for $t \leq d-1$ and $I$ has no generators of degree $\geq d$.
Let $\ell=(\ell_1,\dots,\ell_n)$ be the sequence which satisfies condition (i).
Then 
Propositions \ref{llex} and \ref{lsequence} together with Lemma \ref{hilb}
say that there exists the $d$-linear lexsegment ideal $J' \subset S$ such that
$\ell(J')=\ell$ and $H(J',t)= H(t) $ for all $t \geq d$.
Set 
$$J=J' +I.$$
We will show that this ideal $J$ is a $d$-lexsegment ideal with $H(J,t)=H(t)$ for all $t \in \NN$.

First, we will show $J_{\geq d} =J'$.
Set $A=\{u \in \M_{d-1}:u \in I\}$ and $V=G(J')=\{ u \in \M_d: u \in J'\}$.
What we must prove is
$$\{x_i u :i=1,2,\dots,n \mbox{ and } u \in A\} \subset J'.$$
However, since $J'$ is strongly stable,
it is enough to prove $x_n A \subset J'$.

By condition (i)--(c), we have
\begin{eqnarray}
|A|= H(d-1)\leq \ell_n = \ell_n(J')=|D_n(V)|. \label{LRMQ}
\end{eqnarray}
Since $J'$ is $d$-linear lexsegment,
$D_n(V) \subset \M_{d-1}$ is lexsegment in $K[x_1,\dots,x_n]$.
Since $A \subset \M_{d-1}$ is also lexsegment in $K[x_1,\dots,x_n]$ by the definition of $I$,
(\ref{LRMQ}) says $A \subset D_n(V)$.
Thus we have $x_n A \subset x_n D_n(V) \subset J'$ as desired.

Now, we proved that $J_{\geq d}=J'$ is $d$-linear lexsegment.
Also, since $J'$ is generated in degree $d$,
$J_k=I_k$ is lexsegment for $k<d$.
Thus $J$ is $d$-lexsegment.

It remains to prove that $H(J,t)=H(t)$ for all $t \in \NN$.
It is clear that
$$H(J,t) = \dim_K J_t = \dim_K I_t =H(t) \ \ \mbox{ for } t<d.$$
Also, since $J_{\geq d}=J'$ and $H(J',t)=H(t)$ for all $t \geq d$,
we have
$$H(J,t) = \dim_K J_t = \dim_K J'_t = H(t) \ \ \mbox{ for } t \geq d,$$
as desired.
\end{proof}

Proposition \ref{dbasic} and Theorems \ref{ness} and \ref{suff} immediately imply

\begin{cor} \label{existence}
For any $d$-regular graded ideal $I \subset S$,
there exists the unique $d$-lexsegment ideal, denoted $\Lex^{(d)}(I) \subset S$,
such that $I$ and $\Lex^{(d)}(I)$ have the same Hilbert function.
\end{cor}

We also characterize Hilbert functions of graded ideals $I \subset S$ with $\reg(I)=d$.
We recall the following well known fact.

\begin{lemma}[Bigatti, Hulett and Pardue] \label{bhp}
For any graded ideal $I \subset S$,
one has $\beta_{ij}(I) \leq \beta_{ij}(\Lex(I))$ for all $i$ and $j$.
\end{lemma}

\begin{cor} \label{kurukuru}
If $I \subset S$ is a graded ideal with $\reg(I)=d$ then $\reg(\Lex^{(d)}(I))=d$. 
\end{cor}

\begin{proof}
Suppose that $\reg(\Lex^{(d)}(I)) <d$.
Then Proposition \ref{dbasic} says that  $\Lex^{(d)}(I)=\Lex(I)$.
However, by Lemma \ref{bhp},
one has $d=\reg(I) \leq \reg(\Lex(I))$.
This contradicts the assumption.
\end{proof}

\begin{theorem} \label{equalchara}
Fix a positive integer $d$.
Let $H: \NN \to \NN$ be a numerical function.
Then there exists a nonzero nonunit graded ideal $I \subset S$
such that $\reg(I)=d$ and $H(I,t) = H(t)$ for all $t \in \NN$
if and only if
$H$ satisfies condition $\mathrm{(i)}$ and $\mathrm{(ii)}$
of Theorem \ref{ness} together with the condition 
$\mathrm{(iii)}$ $H(d-1) \up {n-1} < H(d)$.
\end{theorem}

\begin{proof}
(``Only if")
By Corollary \ref{kurukuru},
we may assume that there exists the $d$-lexsegment ideal $I\subset S$
with $\reg(I)=d$ and $H(I,t) =H(t)$ for all $t \in \NN$.
Since $I$ is strongly stable,
$I$ has a generator of degree $d$ by Lemma \ref{ek}.
Thus Lemma \ref{mac} says that $H(I,d-1) \up {n-1} < H(I,d)$.

(``If")
Assume that $H$ satisfies condition (i) and (ii) of Theorem \ref{ness} together with condition (iii).
Then, there exists the $d$-lexsegment ideal $I \subset S$
such that $H(I,t)=H(t)$ for all $t \in \NN$.
We will show $\reg(I)=d$.

Suppose that $\reg(I)<d$.
Then Proposition \ref{dbasic} says that $I$ is lexsegment
and has no generators of degree $\geq d$.
This fact together with Lemma \ref{mac} says that
$H(I,d-1) \up {n-1} =H(I,d)$.
This contradicts condition (iii).
Thus we have $\reg(I)=d$.
\end{proof}

In the rest of this section,
we write two results which follow from Corollary \ref{existence}.
First,
we will give an analogue of Bigatti, Hulett and Pardue's theorem for $d$-lexsegment ideals.

\begin{theorem} \label{maxbetti}
Let $I \subset S$ be a $d$-regular graded ideal.
Then $\beta_{ii+j}(\Lex^{(d)}(I))=\beta_{ii+j}(\Lex(I))$ for all $i$ and $j<d$.
Furthermore, if $\mathrm{char}(K)=0$ then $\beta_{ij}(I) \leq \beta_{ij}(\Lex^{(d)}(I))$
for all $i$ and $j$.
\end{theorem}

\begin{proof}
It follows from \cite[Proposition 2.3]{B} that,
for any strongly stable ideal $J \subset S$,
one has
\begin{eqnarray}
\beta_{ii+k}(J) &=&
\dim_K J_{k} {n-1 \choose i} \label{LRM6}\\
&&-\sum_{q=i}^{n-1} |M_{\leq q}(J,k)| { q-1 \choose i-1} - \sum_{q=i+1}^{n}|M_{\leq q}(J,k-1)| {q-1 \choose i} \nonumber
\end{eqnarray}
for all integers $i$ and $k > 0$,
where $M_{\leq q}(J,k)=\{ u \in \M_k: u \in J \mbox{ and } \max(u) \leq q\}$
for all $1 \leq q \leq n$ and $k > 0$. 
Since $\Lex^{(d)}(I)_k=\Lex(I)_k$ for all $k<d$
the first claim follows from (\ref{LRM6}).

To prove the second claim,
we may assume that $I$ is strongly stable by Lemma \ref{itumono} and Corollary \ref{bs}.
Since $\Lex^{(d)}(I)_k=\Lex(I)_k$ for all $k<d$,
Lemma \ref{bigatti} says that
\begin{eqnarray}
&&|M_{\leq q}(I,k)| \geq |M_{\leq q}(\Lex^{(d)}(I),k)|
\ \ \mbox{ for all $1 \leq q \leq n$ and $k<d$}.  \label{SRMW}
\end{eqnarray}
Since $I$ and $\Lex^{(d)}(I)$ are $d$-regular,
$I_{\geq d}$ and $\Lex^{(d)}(I)_{\geq d}$ are strongly stable ideals generated in degree $d$.
Since they have the same Hilbert function,
we have $\ell(I_{\geq d})=\ell(\Lex^{(d)}(I)_{\geq d})$ by Lemma \ref{equive}.
On the other hand, (\ref{SRM1}) says that
$|M_{\leq q}(I,d)|= \sum_{t=1}^q \ell_t(I_{\geq d})$
and $|M_{\leq q}(\Lex^{(d)}(I),d)|= \sum_{t=1}^q \ell_t(\Lex^{(d)}(I)_{\geq d})$ for all $q$.
Hence we have $|M_{\leq q}(I,d)|= |M_{\leq q}(\Lex^{(d)}(I),d)|$ for all $q$.
By using this fact together with (\ref{LRM6}) and (\ref{SRMW}), we have
$$\beta_{ii+k}(I) \leq \beta_{ii+k}(\Lex^{(d)}(I)) \ \ \mbox{ for all }i \mbox{ and } k \leq d.$$

Since $I$ is $d$-regular,
$\beta_{ii+k}(I) = 0$ for all $i$ and $k >d$.
Thus we have $\beta_{ij}(I) \leq \beta_{ij}(\Lex^{(d)}(I))$ for all $i$ and $j$ as desired.
\end{proof}

It would be interesting to prove the above theorem in arbitrary characteristic.

\begin{exam} \label{exs1}
Let $I=(x_1x_2,x_3x_4) \subset K[x_1,\dots,x_4]$.
Note that $\reg(I) =3$.
Then
$\Lex^{(3)}(I)=(x_1^2,x_1x_2,x_2^2)$,
$\Lex^{(4)}(I)=(x_1^2,x_1x_2,x_1x_3^2,x_2^4)$
and
$\Lex^{(5)}(I)=(x_1^2,x_1x_2,x_1x_3^2,x_1x_3x_4^2,x_2^5,x_2^4x_3).$
(Note that $\Lex^{(6)}(I)=\Lex(I)$.)
The Betti diagrams of these ideals are as follows:\medskip

\begin{center}
$\Lex^{(3)}(I)$: \quad
\begin{tabular}{c|cccccc}
 & 0 & 1 & 2 & 3 \\
\hline
2 & 2 & 1 & - & -\\
3 & 1 & 1 & - & - \\ 
\vspace{-10pt}\\
total & 3 & 2 & 0 & 0\\
\end{tabular}
\quad \quad \quad $\Lex^{(4)}(I)$: \quad
\begin{tabular}{c|cccccc}
 & 0 & 1 & 2 & 3 \\
\hline
2 & 2 & 1 & - & -\\
3 & 1 & 2 & 1 & -\\
4 & 1 & 1 & - & -\\
\vspace{-10pt}\\
total & 4 & 4 & 1 & 0\\
\end{tabular}
\\
\medskip
\end{center}
\hspace{26pt} $\Lex^{(5)}(I)$: \quad
\begin{tabular}{c|ccccccc}
 & 0 & 1 & 2 & 3 \\
\hline
2 & 2 & 1 & - & -\\
3 & 1 & 2 & 1 & -\\
4 & 1 & 3 & 3 & 1\\
5 & 2 & 3 & 1 & -\\
\vspace{-10pt}\\
total & 6 & 9 & 5 & 1\\
\end{tabular}
\end{exam}

Let $H: \mathbb{N} \to \mathbb{N}$ be a numerical function.
We say that $H$ is \textit{admissible} if there exists a graded ideal $I \subset S$ such that $H(S/I,t)=H(t)$ for all $t \in \NN$.
For an admissible numerical function $H$,
set 
$$\A_H=\{\mathrm{depth}(S/I): I \mbox{ is a graded ideal with } H(S/I,t)=H(t)
\mbox{ for all } t \in \NN\},$$
where $\mathrm{depth}(S/I)$ is the depth of $S/I$, and
$$\R_H=\{\reg(I): I \mbox{ is a graded ideal with } H(S/I,t)=H(t)
\mbox{ for all } t \in \NN\}.$$
In \cite{MH},
it was proved that $\A_H$ is a set of integers of the form
$\A_H=\{a,a+1,\dots,b-1,b\}$ for some $0 \leq a\leq b$.
We will give an analogue of this fact for $\R_H$.

\begin{theorem} \label{regdence}
Let $H$ be an admissible function,
$a=\min \R_H$ and $b=\max \R_H$.
For any integer $a \leq r \leq b$,
one has $r \in \R_H$.
\end{theorem}

\begin{proof}
Let $I$ be a graded ideal with $\reg(I)=a$ and $H(S/I,t)=H(t)$
for all $t \in \NN$.
Fix an integer $a \leq r \leq b$.
Since $I$ is $r$-regular,
there exists the $r$-lexsegment ideal $\Lex^{(r)}(I)$
with the same Hilbert function as $I$.
We will show $\reg(\Lex^{(r)}(I))=r$.

Suppose that $\reg(\Lex^{(r)}(I))<r$.
Then, by Proposition \ref{dbasic},
we have $\Lex^{(r)}(I)=\Lex(I)$.
However, by Lemma \ref{bhp},
$\reg(\Lex(I))=b$.
Since we assume $r \leq b$, this is a contradiction.
Thus $\reg(\Lex^{(r)}(I))=r$.
\end{proof}

\begin{exam}
Let $I=(x_1x_2,x_3x_4) \subset S=K[x_1,\dots,x_4]$,
and let $H$ be the Hilbert function of $S/I$.
Then, by the computations in Example \ref{exs1}, we have
$\A_H=\{0,1,2\}$ and $\R_H=\{3,4,5,6\}$.
\end{exam}

\section{Squarefree $d$-lexsegment ideals}

In this section,
we define $d$-linear squarefree lexsegment ideals
and squarefree $d$-lexsegment ideals,
and study their properties.

We say that a set of squarefree monomials $V \subset S$ is \textit{squarefree lexsegment}
if, for all squarefree monomial ideals $u \in V$
and $v >_\lex u$ with $\deg(v)=\deg(u)$,
it follows that $v \in V$.
A squarefree monomial ideal $I \subset S$ is said to be \textit{squarefree lexsegment}
if the set of squarefree monomials in $I$ is squarefree lexsegment.

Let $1 \leq d \leq n$ be an integer, $\M \sqcomp d$ the set of all squarefree monomials in $S$
of degree $d$,
$V \subset \M \sqcomp d$ a squarefree strongly stable set of squarefree monomials
and $I$ the squarefree monomial ideal generated by $V$.
Recall that, in \S 3,
we define $\ell^*(I)=(\ell^*_1(I),\dots,\ell^*_{n-d+1}(I))$
by $\ell^*_k(I)=|\{u \in G(I): \max(u)=k+d-1\}|$
and
$$D_k(V)=\{(u / x_k) \in \M \sqcomp {d-1}: u \in V \mbox{ and } \max(u)=k\} \ \ \mbox{ for }k \geq 1.$$
Set $\ell^*(V)=\ell^*(I)$.
Then it is easy to see that
$$\ell^*_k(V)=|D_{k+d-1}(V)| \ \ \mbox{ for }k=1,2,\dots,n-d+1.$$
The following fact follows from Lemma \ref{sqek} in the same way as Lemma
\ref{equive}.

\begin{lemma} \label{sqequive}
Let $I \subset S$ and $J \subset S$ be squarefree strongly stable ideals
generated in degree $d$.
Then the following conditions are equivalent.
\begin{itemize}
\item[(i)] $I$ and $J$ have the same Hilbert function;
\item[(ii)] $I$ and $J$ have the same graded Betti numbers;
\item[(iii)] $\ell^*(I)= \ell^*(J)$.
\end{itemize}
\end{lemma}

A subset
$V \subset \M \sqcomp d$
is said to be \textit{$d$-linear squarefree lexsegment}
if $V$ is squarefree strongly stable and each $D_k(V)$ is squarefree lexsegment in $K[x_1,\dots,x_{k-1}]$
for $k=d,d+1,\dots,n$.
Also,
a squarefree monomial ideal $I \subset S$
is said to be \textit{$d$-linear squarefree lexsegment} if $I$
is generated by a $d$-linear squarefree lexsegment subset $V \subset \M \sqcomp d$.
The squarefree operation gives a nice relation between $d$-linear lexsegment
ideals and $d$-linear squarefree lexsegment ideals.
Indeed, we have

\begin{lemma} \label{corres}
Let $d$ be a positive integer.
A monomial ideal $I \subset K[x_1,\dots,x_n]$ is $d$-linear lexsegment
if and only if
$\Phi(I) \subset K[x_1,\dots,x_{n+d-1}]$
is $d$-linear squarefree lexsegment.
\end{lemma}

\begin{proof}
By Lemma \ref{operation},
$I \subset K[x_1,\dots,x_n]$ is a strongly stable ideal generated in degree $d$
if and only if $\Phi(I) \subset K[x_1,\dots,x_{n+d-1}]$ is a squarefree
strongly stable ideal generated in degree $d$.
Also, by the definition of the squarefree operation, we have
$$D_{k+d-1}(G(\Phi(I)))=\{ \Phi(u): u \in D_k(G(I))\} \ \ \mbox{ for }
k=1,2,\dots,n.$$
Thus what we must prove is $D_k(G(I))$ is lexsegment if and only if
$\{ \Phi(u): u \in D_k(G(I))\}$ is squarefree lexsegment.
However, this statement follows from the fact that,
for any monomials $u,v \in \M_d$,
one has $u >_\lex v$ if and only if $\Phi(u) >_\lex \Phi(v)$.
\end{proof}

\begin{lemma} \label{lsqlex}
Fix an integer $1 \leq d \leq n$.
For any squarefree strongly stable ideal $I \subset S$
generated in degree $d$,
there exists a $d$-linear squarefree lexsegment ideal $J \subset S$
such that $\ell^* (I)=\ell^* (J)$.
\end{lemma}

\begin{proof}
By Lemma \ref{operation},
$\Phi^{-1}(I) \subset K[x_1,\dots,x_{n-d+1}]$ is a strongly stable ideal generated in degree $d$.
Then Proposition \ref{llex} says that there exists the $d$-linear lexsegment ideal
$J' \subset K[x_1,\dots,x_{n-d+1}]$ with $\ell(J')=\ell(\Phi^{-1}(I))$.
Set $J=\Phi(J') \subset S$.
Then $J$ is $d$-linear squarefree lexsegment by Lemma \ref{corres}.
Since $\ell^*(I)=\ell(\Phi^{-1}(I))$ and $\ell(J')=\ell^*(\Phi(J'))=\ell^*(J)$,
we have $\ell^*(I)=\ell^*(J)$ as required.
\end{proof}

\begin{definition}
Let $1 \leq d \leq n$ be an integer.
For an integer $1 \leq k \leq n$ and a squarefree monomial ideal $I \subset S$,
we say that $I_k$ is \textit{squarefree lexsegment} if the set of squarefree monomials in $I_k$ is squarefree lexsegment,
and write
$I \sqcompleq k$ for the ideal generated by all squarefree monomials
in $I$ of degree $\geq k$.
A squarefree monomial ideal $I \subset S$ is said to be
\textit{squarefree $d$-lexsegment}
if $I \sqcompleq d$ is $d$-linear squarefree lexsegment and
$I_k$ is squarefree lexsegment for all $1 \leq k <d$.
\end{definition}

The next lemma can be proved in the same way as Proposition \ref{dbasic}.

\begin{proposition} \label{sqdbasic}
Let $1 \leq d \leq n$ be an integer.
\begin{itemize}
\item[(i)] If $I \subset S$ is squarefree $d$-lexsegment then $I$ is 
squarefree strongly stable and $d$-regular;
\item[(ii)] If $I$ and $J$ are squarefree $d$-lexsegment and
have the same Hilbert function then $I=J$;
\item[(iii)] If $I$ is a squarefree $d$-lexsegment ideal with $\reg(I)<d$
then $I$ is squarefree lexsegment.
\end{itemize}
\end{proposition}

\begin{proof}
(i) It is clear that $I$ is squarefree strongly stable.
Then $I$ is $d$-regular by Lemma \ref{sqek} since $I$ has no generators of degree $> d$.

(ii) 
It is enough to show that the set of squarefree monomials in $I$ is equal to that in $J$.
Since $I_k$ and $J_k$ are squarefree lexsegment for $k<d$ and $I$ and $J$
have the same Hilbert function,
the set of squarefree monomials in $I_k$ is equal to that in $J_k$ for $k<d$.
On the other hand, $I \sqcomp {\geq d}$ and $J \sqcomp{ \geq d}$ are $d$-linear squarefree
lexsegment and have the same Hilbert function.
Since $d$-linear squarefree lexsegment ideals are uniquely determined from their $\ell^*$-sequence,
we have $I \sqcompleq d = J \sqcompleq d$ by Lemma \ref{sqequive}.
Hence we have $I=J$.

(iii) Since $I$ is squarefree strongly stable, $\reg(I)<d$ implies
that $I$ has no generators of degree $\geq d$ by Lemma \ref{sqek}.
Then $I = I_{\leq d-1}$ is squarefree lexsegment.
\end{proof}

A simplicial complex $\Gamma$ satisfying $\reg(I_\Gamma)\leq d$ is called a $(d-1)$-Leray
simplicial complex.
In 1973 Eckhoff proposed a conjecture about the characterization of $f$-vectors of Leray simplicial complexes.
The necessity of the conjecture was proved by Kalai \cite{K10},
and the sufficiency of the conjecture was proved by
Eckhoff and Kalai \cite{K11} independently.
(Eckhoff did not publish the proof.)
About the precise conditions of Eckhoff and Kalai's result, see \cite{K10}.
The characterization of Hilbert functions of squarefree monomial ideals
which are $d$-regular follows from their result.
We will prove that,
for any squarefree monomial ideal $I \subset S$ with $\reg(I) \leq d$,
there exists the squarefree $d$-lexsegment ideal with the same Hilbert function as $I$.
Of course this fact follows from Eckhoff and Kalai's result by checking that
if $\Gamma$ is a simplicial complex which was used in the proof of the sufficiency of their result \cite{K11}
then $I_\Gamma$ is squarefree $d$-lexsegment.
Here, we will give a different proof which is similar to the proof
of Theorem \ref{suff}.
The next fact follows from the Kruskal--Katona Theorem
(see \cite[\S 4]{AHH}).

\begin{lemma}[Kruskal--Katona] \label{kk}
For any squarefree monomial ideal $I \subset S$,
there exists the unique squarefree lexsegment ideal,
denoted $\SqLex(I) \subset S$,
such that $I$ and $\SqLex(I)$ have the same Hilbert function.
\end{lemma}

\begin{lemma}[{\cite[Theorem 3.9]{SQlex}}] \label{sqbig}
Let $L \subset \M \sqcomp d$ be squarefree lexsegment,
and let $V \subset \M \sqcomp d$ be a strongly stable set of 
squarefree monomials with $|V|=|L|$.
Then
$$|M_{\leq k}(V)| \geq |M_{\leq k}(L)|\ \ \mbox{ for all }k.$$
\end{lemma}

\begin{proposition}\label{sqexist}
Fix an integer $1 \leq d \leq n$.
Let $I \subset S$ be a squarefree monomial ideal with $\reg(I) \leq d$.
Then there exists the unique squarefree $d$-lexsegment ideal,
denoted $\SqLex^{(d)}(I) \subset S$,
such that $I$ and $\SqLex^{(d)}(I)$ have the same Hilbert function. 
\end{proposition}

\begin{proof}
We may assume that $I$ is squarefree strongly stable by Lemma \ref{shifting}.
Then Lemma \ref{kk} says that
there exists the squarefree lexsegment ideal $I' \subset S$
such that $\dim_K I_k = \dim_K I'_k$ for $k<d$
and $I'$ has no generators of degree $\geq d$.
Also, since $I$ is $d$-regular,
$I \sqcompleq d$ is a squarefree strongly stable ideal generated in degree $d$ by Lemma \ref{sqek}.
Then Lemma \ref{lsqlex} says that
there exists the $d$-linear squarefree lexsegment ideal $J' \subset S$ such that
$\ell^*(J')=\ell^*(I \sqcompleq d)$.
Set
$$J=J' + I'.$$ 
We will show that $J$ is a squarefree $d$-lexsegment ideal which has the same Hilbert function as $I$.
First, we will show $J \sqcompleq d =J'$.
Let $A=\{u \in \M \sqcomp {d-1}:u\in I'\}$ and
$V=\{u \in \M \sqcomp d : u \in J'\}$.
What we must prove is $\{x_iu\in \M \sqcomp d: u \in A \mbox{ and }i=1,2,\dots,n\} \subset V$.
However, since $A$ and $V$ are squarefree strongly stable,
it is enough to prove that $x_n M_{\leq n-1}(A) \subset V$.

Let $W_k=\{u \in \M \sqcomp k: u \in I\}$ for $k \geq 0$.
It is clear that $x_n M_{\leq n-1}(W_{d-1}) \subset W_d$.
Hence $M_{\leq n-1}(W_{d-1}) \subset D_n(W_d)$,
where $D_n(W_d)=\{u/x_n: u \in W_d,\ \max(u)=n\}$.
Since $A$ is a squarefree lexsegment subset of $\M \sqcomp d$ with $|A|=\dim I_{d-1} = |W_{d-1}|$,
Lemma \ref{sqbig} says
$$|M_{\leq n-1}(A)| \leq |M_{\leq n-1}(W_{d-1})| \leq  |D_n (W_d)|.$$
On the other hand,
since $\ell^*(V)=\ell^*(J')=\ell^*(I \sqcompleq d) =\ell^*(W_d)$,
it follows that
$$|D_n(W_d)|=\ell_n^*(I \sqcompleq d)=\ell_n^*(V)=|D_n(V)|.$$
Thus we have
$|M_{\leq n-1}(A)| \leq |D_n(V)|$.
Since $J'$ is $d$-linear squarefree lexsegment,
$D_n(V)$ is squarefree lexsegment in $K[x_1,\dots,x_{n-1}]$.
Since $A$ is squarefree lexsegment,
$M_{\leq n-1}(A)$ is also squarefree lexsegment in $K[x_1,\dots,x_{n-1}]$.
Thus we have $M_{\leq n-1}(A) \subset D_n(V)$.
Hence we have $x_n  M_{\leq n-1}(A)  \subset V$.
Thus $J \sqcompleq d =J'$.

Then,
by the construction,
$J \sqcompleq d =J'$ is $d$-linear squarefree lexsegment
and $J_k=I_k'$ is squarefree lexsegment for $k<d$.
Thus $J$ is squarefree $d$-lexsegment.
To prove that $I$ and $J$ have the same Hilbert function,
it is enough to show that the number of squarefree monomials in $I$ of degree $k$
is equal to that in $J$ for all $1 \leq k \leq n$.
Since $\dim I_k = \dim_K I'_k=\dim_K J_k$ for $k<d$,
the number of squarefree monomials in $J_k$ is equal to that in $I_k$ for $ k<d$.
On the other hand,
since $J \sqcompleq d =J'$ and $\ell^*( I \sqcompleq d)=\ell^* (J')$,
Lemma \ref{sqequive} says that the number of squarefree monomial ideals in $J$ of degree $k$
is also equal to that in $I$ for $k \geq d$.
Thus $I$ and $J$ have the same Hilbert function.
\end{proof}

Next, we prove an analogue of Theorem \ref{regdence}.
A numerical function $H:  \NN \to \NN$ is said to be a \textit{squarefree function}
if there exists a squarefree monomial ideal $I \subset S$ such that $H(S/I,t)=H(t)$ for all $t \in \NN$.
Let $\mathrm{Sq}(S)$ be the set of all squarefree monomial ideals in $S$
and define
$$\mathcal{SR}_H=\{\reg(I): I \in \mathrm{Sq}(S),\ H(S/I,t)=H(t) \mbox{ for all } t \in \NN\}.$$
The following fact was first proved in \cite{AHHshifting} 
for base fields of characteristic $0$, and
a proof in any characteristic was later given in \cite{MHshifting}.

\begin{lemma} \label{sqbhp}
If $I \subset S$ is a squarefree monomial ideal
then $\beta_{ij}(I) \leq \beta_{ij}(\SqLex(I))$ for all $i$ and $j$.
\end{lemma}



\begin{theorem}
\label{monowasure}
Let $1 \leq d \leq n$ be an integer and $H: \NN \to \NN$
a squarefree function.
Set $a = \min \mathcal{SR}_H$ and $b= \max \mathcal{SR}_H$.
Then, for any integer $a \leq r \leq b$,
one has $r \in \mathcal{SR}_H$.
\end{theorem}

\begin{proof}
Let $I \subset S$ be a squarefree monomial ideal with $\reg(I)=a$ and $H(S/I,t)=H(t)$ for all $t \in \NN$.
Fix an integer $a \leq r \leq b$.
We will show that $\reg(\SqLex^{(r)}(I))=r$.

Suppose that $\reg(\SqLex^{(r)}(I))<r$.
Then $\SqLex^{(r)}(I)=\SqLex(I)$ by Lemma \ref{sqdbasic}.
However, by Lemma \ref{sqbhp}, we have $r>\reg(\SqLex(I))=b$.
Since we assume $r\leq b$, this is a contradiction.
Thus $\reg(\SqLex^{(r)}(I))=r$.
\end{proof}

Finally, we will give an analogue of Theorem \ref{maxbetti}.
We first recall some results on the squarefree operation which were used in \cite{AHHshifting}.

It is well known that any squarefree monomial ideal $I \subset S$ satisfies
$\beta_{ij}(I) =0$ for all $i$ and $j >n$.
On the other hand,
if $I$ is a strongly stable ideal satisfying $\beta_{ij}(I)=0$
for all $i$ and $j>n$,
then the Eliahou--Kervaire formula says
$\max(u) + \deg(u) -1 \leq n$ for any $u \in G(I)$.
Also, if a monomial $u \in S$ satisfies $\max(u) + \deg(u) -1 \leq n$
then $\max(\Phi(u)) \leq n$.
Thus, if $I \subset S$ is a strongly stable ideal satisfying $\beta_{ij}(I)=0$
for all $i$ and $j>n$, then we write $\tilde \Phi(I) \subset S$ for the ideal generated by $\{ \Phi(u): u\in G(I)\}$.

\begin{lemma}[Aramova--Herzog--Hibi \cite{AHHshifting}] \label{selfoperator}
If $I \subset S$ is a strongly stable ideal satisfying $\beta_{ij}(I)=0$
for all $i$ and $j>n$,
then $\tilde \Phi(I)\subset S$ is a squarefree strongly stable ideal
and has the same graded Betti numbers as $I$.
\end{lemma}

Let $I \subset S$ be a squarefree monomial ideal.
If $\mathrm{char}(K)=0$ then
Lemmas \ref{itumono} and \ref{bcp} says that $\gin(I)$ is a strongly stable ideal
satisfying $\beta_{ij}(\gin(I))=0$ for all $i$ and $j>n$.
Thus we can define $\tilde \Phi(\gin(I)) \subset S$.
The next lemma immediately follows from Lemma \ref{selfoperator}.

\begin{lemma} \label{tsuka}
Assume $\mathrm{char}(K)=0$.
If $I \subset S$ is a squarefree monomial ideal then
$\tilde \Phi(\gin(I))\subset S$ is a squarefree strongly stable ideal
satisfying 
$\beta_{ij}(\tilde \Phi(\gin(I))) =\beta_{ij}(\gin(I))$ for all $i$
and $j$.
\end{lemma}

Note that
the operation $I \to \tilde \Phi(\gin(I))$ was considered by Kalai
\cite{K2} to define
symmetric algebraic shifting.
We do not give a definition of symmetric algebraic shifting here.
See, e.g., \cite{H} or \cite{K}.

\begin{theorem} \label{sqmaxbetti}
Let $1 \leq d \leq n$ be an integer and $I \subset S$
a squarefree monomial ideal with $\reg(I) \leq d$.
Then $\beta_{ii+j}(\SqLex^{(d)}(I))=\beta_{ii+j}(\SqLex(I))$
for all $i$ and $j<d$.
Furthermore, if $\mathrm{char}(K)=0$ then 
$\beta_{ij}(I) \leq \beta_{ij}(\SqLex^{(d)}(I))$ for all $i$ and $j$.
\end{theorem}

\begin{proof}
It follows from \cite[Theorem 4.4]{SQlex} that, for any squarefree strongly stable ideal $J\subset S$, we have
\begin{eqnarray}
\ \ \beta_{ii+k}(J) &=& \dim_K J_k { n-k \choose i}  \label{parasol} \\
&&\nonumber -\bigg\{\sum_{t=k}^{n-1}|M^*_{\leq t}(J,k)|{t-k \choose i-1} +\sum_{t=k}^n|M^*_{\leq t-1}(J,k-1)|{t-k \choose i}\bigg\}
\end{eqnarray}
for all $i$ and $k$,
where $M^*_{\leq t}(J,k)=\{ u \in J: u \in \M \sqcomp k\mbox{ and } \max(u) \leq t\}$.
Since $\SqLex^{(d)}(I)_k=\SqLex(I)_k$ for all $k<d$,
the first claim follows from  (\ref{parasol}).

To prove the second claim,
we may assume that $I$ is squarefree strongly stable
by Lemmas \ref{bcp} and \ref{tsuka}.
Since $\SqLex^{(d)}(I)_k=\SqLex(I)_k$ for all $k<d$,
Lemma \ref{sqbig} says that
\begin{eqnarray}
M^*_{\leq t}(I,k) \geq  M^*_{\leq t}(\SqLex^{(d)}(I),k)
\ \ \mbox{ for all } t \mbox{ and }  k < d. \label{LRMA}
\end{eqnarray}
On the other hand,
since $I$ and $\SqLex^{(d)}(I)$ are $d$-regular,
$I \sqcompleq d$ and $\SqLex^{(d)}(I) \sqcompleq d$
are squarefree strongly stable ideals generated in degree $d$
by Lemma \ref{sqek}.
Notice that $I \sqcompleq d$ and $\SqLex^{(d)}(I) \sqcompleq d$ have the same Hilbert function.
Then, in the same way as the proof of Theorem \ref{maxbetti},
Lemma \ref{sqequive} together with the definition of $\ell^*$-sequences imply
\begin{eqnarray}
&&M^*_{\leq t}(I,d) = \sum_{k=1}^t \ell^*_k(I \sqcompleq d)
= \sum_{k=1}^t \ell^*_k(\SqLex^{(d)}(I) \sqcompleq d)= M^*_{\leq t}(\SqLex^{(d)}(I),d)
\label{LRMB}
\end{eqnarray}
for all $t$.
Then (\ref{parasol}), (\ref{LRMA}) and (\ref{LRMB}) say that
$\beta_{ii+k}(I) \leq \beta_{ii+k}(\SqLex^{(d)}(I))$
for all $i$ and $ k \leq d$.
Since $\reg(I) \leq d$ we have $\beta_{ii+k}(I) =0$ for all $i$ and $k >d$.
Thus the claim follows.
\end{proof}

\begin{exam}
Let $I=(x_1x_3x_5,x_1x_3x_6,x_1x_4x_6,x_2x_4x_6) \subset S=K[x_1,\dots,x_6]$.
Note that $\reg(I)=3$.
Then 
$$\SqLex^{(3)}(I)=
(x_1x_2x_3,x_1x_2x_4,x_1x_3x_4,x_2x_3x_4),$$
$$ \SqLex^{(4)}(I)=
(x_1x_2x_3,x_1x_2x_4,x_1x_2x_5,x_1x_2x_6,x_1x_3x_4x_5,x_1x_3x_4x_6,
x_2x_3x_4x_5)$$
and $\SqLex(I)$ is
\begin{eqnarray*}
(x_1x_2x_3,x_1x_2x_4,x_1x_2x_5,x_1x_2x_6,x_1x_3x_4x_5,x_1x_3x_4x_6,
x_1x_3x_5x_6,x_2x_3x_4x_5x_6).
\end{eqnarray*}
Betti diagrams  are following:

\begin{center}
$\SqLex^{(3)}(I)$: \quad
\begin{tabular}{c|ccccc}
 & 0 & 1 & 2 & 3 \\
\hline
3 & 4 & 3 & - & -\\
\vspace{-10pt}\\ 
total & 4\ & 3\ & 0 & 0\\
\end{tabular}
\quad \quad 
$\SqLex^{(4)}(I)$: \quad
\begin{tabular}{c|ccccc}
 & 0 & 1 & 2 & 3 \\
\hline
3 & 4 & 6 & 4 & 1\\
4 & 3 & 4 & 1 & -\\
\vspace{-10pt}\\ 
total & 7\ & 10 & 5 & 1\\
\end{tabular}
\end{center}

\hspace{8pt}$\SqLex(I)$: \quad \quad
\begin{tabular}{c|ccccc}
 & 0 & 1 & 2 & 3 \\
\hline
3 & 4 & 6 & 4 & 1\\
4 & 3 & 5 & 2 & -\\
5 & 1 & 1 & - & -\\
\vspace{-10pt}\\ 
total & 8 & 12 & 6 & 1\\
\end{tabular}
\end{exam}

\begin{rem}
Since $\reg(I_\Gamma)=n-\mathrm{depth}(S/I_{\Gamma^*})$,
the characterization of $f$-vectors of simplicial complexes $\Gamma$ with $\reg(I_\Gamma) \leq n-s$
is equivalent to that of simplicial complexes $\Gamma$ with $\mathrm{depth}(S/I_\Gamma) \geq s$,
which was presented by Bj\"orner \cite[Corollary 3.2]{Bj}.
Since squarefree $d$-lexsegment ideals characterize the Hilbert functions of squarefree monomial ideals $I$ of $S$ with $\reg(I) \leq d$,
by using the Alexander duality,
Proposition \ref{sqexist}
is essentially equivalent to Bj\"orner's result.
Also, by considering the Alexander dual of a simplicial complex whose Stanley--Reisner ideal is
squarefree $d$-lexsegment,
we can prove analogues of Propositions \ref{sqdbasic}, \ref{sqexist} and Theorem \ref{monowasure} for depth.
Note that an analogue of Theorem \ref{sqmaxbetti} for depth will be proved in the next section
(see Corollary \ref{wasure}).
\end{rem}

\section{A generalization of Bigatti, Hulett and Pardue's theorem}

Throughout this section,
we assume that the base field $K$ is a field of characteristic $0$.
In this section,
we will consider strongly stable ideals having the maximal
graded Betti numbers
among graded ideals which
have the same Hilbert function and
which accept certain restrictions
on their graded Betti numbers.

\begin{definition}
Set $\set=\{0,1,2,\dots,n-1\}$.
A finite subset $A \subset \area$ is said to be an \textit{extremal area} if
$(i,j) \in A$ implies $(i',j') \in A$ for all $(i',j') \in \area$
with $i' \leq i$ and $j' \leq j$.
For any $(i,j) \in \area$,
let $\langle (i,j) \rangle$ denote the minimal extremal area which contains $(i,j)$,
in other words, 
$$\langle (i,j) \rangle=\{ (i',j') \in \area:i' \leq i \mbox{ and }
j' \leq j\}.$$
It is easily verified that any extremal area $A \subset \area$
has the unique representation
$$A= \bigcup_{k=1}^t \langle (i_k,j_k) \rangle$$
satisfying $t \leq n$, $i_1 < \cdots < i_t$ and $j_1 > \cdots >j_t$.
This unique representation will be called the \textit{standard representation} of $A$
and the elements $(i_1,j_1),\dots,(i_t,j_t)$ will be called
\textit{extremal points of $A$}.
(See the picture bellow.)
\end{definition}
\pagebreak

{
\begin{center}
An extremal area $A=\bigcup_{k=1}^3 \langle (i_k,j_k) \rangle$ \medskip\\
\unitlength 0.1in
\begin{picture}( 38.7500, 18.0000)(  7.7500,-20.0000)
\put(15.0000,-7.0000){\makebox(0,0){$\bullet$}}%
\put(15.0000,-9.0000){\makebox(0,0){$\bullet$}}%
\put(15.0000,-11.0000){\makebox(0,0){$\vdots$}}%
\put(15.0000,-14.0000){\makebox(0,0){$\bullet$}}%
\put(15.0000,-16.0000){\makebox(0,0){$\bullet$}}%
\put(15.0000,-18.0000){\makebox(0,0){$\bullet$}}%
\put(18.0000,-7.0000){\makebox(0,0){$\cdots$}}%
\put(18.0000,-9.0000){\makebox(0,0){$\cdots$}}%
\put(18.0000,-11.5000){\makebox(0,0){$\cdots$}}%
\put(18.0000,-14.0000){\makebox(0,0){$\cdots$}}%
\put(18.0000,-16.0000){\makebox(0,0){$\cdots$}}%
\put(18.0000,-18.0000){\makebox(0,0){$\cdots$}}%
\put(21.0000,-18.0000){\makebox(0,0){$*$}}%
\put(21.0000,-16.0000){\makebox(0,0){$\bullet$}}%
\put(21.0000,-14.0000){\makebox(0,0){$\bullet$}}%
\put(21.0000,-11.0000){\makebox(0,0){$\vdots$}}%
\put(21.0000,-9.0000){\makebox(0,0){$\bullet$}}%
\put(21.0000,-7.0000){\makebox(0,0){$\bullet$}}%
\put(23.5000,-14.0000){\makebox(0,0){$\bullet$}}%
\put(23.5000,-11.0000){\makebox(0,0){$\vdots$}}%
\put(23.5000,-9.0000){\makebox(0,0){$\bullet$}}%
\put(23.5000,-7.0000){\makebox(0,0){$\bullet$}}%
\put(26.0000,-14.0000){\makebox(0,0){$\bullet$}}%
\put(26.0000,-11.0000){\makebox(0,0){$\vdots$}}%
\put(26.0000,-9.0000){\makebox(0,0){$\bullet$}}%
\put(26.0000,-7.0000){\makebox(0,0){$\bullet$}}%
%
\special{pn 8}%
\special{pa 3900 400}%
\special{pa 850 400}%
\special{fp}%
\put(12.5000,-7.0000){\makebox(0,0){$\bullet$}}%
\put(12.5000,-9.0000){\makebox(0,0){$\bullet$}}%
\put(12.5000,-11.0000){\makebox(0,0){$\vdots$}}%
\put(12.5000,-14.0000){\makebox(0,0){$\bullet$}}%
\put(12.5000,-16.0000){\makebox(0,0){$\bullet$}}%
\put(12.5000,-18.0000){\makebox(0,0){$\bullet$}}%
\put(12.5000,-5.0000){\makebox(0,0){$\bullet$}}%
\put(15.0000,-5.0000){\makebox(0,0){$\bullet$}}%
\put(18.0000,-5.0000){\makebox(0,0){$\cdots$}}%
\put(21.0000,-5.0000){\makebox(0,0){$\bullet$}}%
\put(23.5000,-5.0000){\makebox(0,0){$\bullet$}}%
\put(26.0000,-5.0000){\makebox(0,0){$\bullet$}}%
%
\special{pn 8}%
\special{pa 1150 200}%
\special{pa 1150 2000}%
\special{fp}%
\put(28.5000,-11.0000){\makebox(0,0){$\vdots$}}%
\put(28.5000,-9.0000){\makebox(0,0){$\bullet$}}%
\put(28.5000,-7.0000){\makebox(0,0){$\bullet$}}%
\put(28.5000,-5.0000){\makebox(0,0){$\bullet$}}%
\put(28.5000,-14.0000){\makebox(0,0){$*$}}%
\put(31.0000,-5.0000){\makebox(0,0){$\bullet$}}%
\put(31.0000,-7.0000){\makebox(0,0){$\bullet$}}%
\put(31.0000,-9.0000){\makebox(0,0){$\bullet$}}%
\put(34.0000,-5.0000){\makebox(0,0){$\cdots$}}%
\put(34.0000,-7.0000){\makebox(0,0){$\cdots$}}%
\put(34.0000,-9.0000){\makebox(0,0){$\cdots$}}%
\put(37.0000,-9.0000){\makebox(0,0){$*$}}%
\put(37.0000,-7.0000){\makebox(0,0){$\bullet$}}%
\put(37.0000,-5.0000){\makebox(0,0){$\bullet$}}%
\put(10.0000,-5.0000){\makebox(0,0){$1$}}%
\put(10.0000,-9.0000){\makebox(0,0){$j_3$}}%
\put(10.0000,-14.0000){\makebox(0,0){$j_2$}}%
\put(10.0000,-18.0000){\makebox(0,0){$j_1$}}%
\put(12.5000,-3.0000){\makebox(0,0){$0$}}%
\put(15.0000,-3.0000){\makebox(0,0){$1$}}%
\put(21.0000,-3.0000){\makebox(0,0){$i_1$}}%
\put(28.5000,-3.0000){\makebox(0,0){$i_2$}}%
\put(37.0000,-3.0000){\makebox(0,0){$i_3$}}%
%
\special{pn 8}%
\special{pa 3850 400}%
\special{pa 3850 1000}%
\special{pa 3000 1000}%
\special{pa 3000 1500}%
\special{pa 2250 1500}%
\special{pa 2250 1900}%
\special{pa 1150 1900}%
\special{fp}%
\end{picture}%
\end{center}
(Each $\bullet$ denotes an element in $A$ and $*$ denotes an extremal point of $A$.)
\medskip
}

Let $A \subset \area$ be an extremal area.
We say that a graded ideal $I \subset S$ \textit{admits
the extremal area $A$} if $\beta_{ii+j}(I)=0$ for all $(i,j) \not \in A$.
For a numerical function $H:\NN \to \NN$,
define the set of graded ideals $\mathcal{L}_H^A$ by
$$\mathcal{L}_H^A=\{ I \subset S:I \mbox{ admits } A
\mbox{ and } H(S/I,t)=H(t) \mbox{ for all } t \in \NN\}.$$

Since $I$ and $\gin(I)$ have the same extremal Betti numbers,
if $I$ admits $A$ then $\gin(I)$ also admits $A$.
One may expect the existence of a strongly stable ideal $L \in \mathcal{L}_H^A$ satisfying $\beta_{ij}(L) \geq \beta_{ij}(J)$
for all $i,j$ and $J \in \mathcal{L}_H^A$.
Unfortunately,
this statement is false for some extremal areas.

\begin{exam} \label{counter}
Let $S=K[x_1,\dots,x_5]$,
$I=(x_1^2,x_1x_2,x_1x_3,x_1x_4,x_2^2,x_2x_3^3,x_3^4) \subset S$
and $H$ the Hilbert function of $S/I$.
Set $A=\langle (2,4) \rangle \cup \langle (4,2) \rangle$.
Then the Betti diagram of $I$ is the following.
(The line in the diagram corresponds to the extremal area $A$.)

\begin{center}
\unitlength 0.1in
\begin{picture}( 14.8500, 11.0000)( 13.6500,-17.0000)
\put(17.0000,-9.0000){\makebox(0,0){-}}%
\put(17.0000,-11.0000){\makebox(0,0){$5$}}%
\put(17.0000,-13.0000){\makebox(0,0){-}}%
\put(17.0000,-15.0000){\makebox(0,0){$2$}}%
\put(19.5000,-15.0000){\makebox(0,0){$4$}}%
\put(19.5000,-13.0000){\makebox(0,0){-}}%
\put(19.5000,-11.0000){\makebox(0,0){$7$}}%
\put(19.5000,-9.0000){\makebox(0,0){-}}%
\put(19.5000,-7.0000){\makebox(0,0){$1$}}%
\put(22.0000,-7.0000){\makebox(0,0){$2$}}%
\put(22.0000,-9.0000){\makebox(0,0){-}}%
\put(22.0000,-11.0000){\makebox(0,0){$4$}}%
\put(22.0000,-13.0000){\makebox(0,0){-}}%
\put(22.0000,-15.0000){\makebox(0,0){$2$}}%
\put(24.5000,-15.0000){\makebox(0,0){-}}%
\put(24.5000,-13.0000){\makebox(0,0){-}}%
\put(24.5000,-11.0000){\makebox(0,0){$1$}}%
\put(24.5000,-9.0000){\makebox(0,0){-}}%
\put(24.5000,-7.0000){\makebox(0,0){$3$}}%
\put(27.0000,-7.0000){\makebox(0,0){$4$}}%
\put(27.0000,-9.0000){\makebox(0,0){-}}%
\put(27.0000,-11.0000){\makebox(0,0){-}}%
\put(27.0000,-13.0000){\makebox(0,0){-}}%
\put(27.0000,-15.0000){\makebox(0,0){-}}%
\put(17.0000,-7.0000){\makebox(0,0){$0$}}%
\put(15.0000,-9.0000){\makebox(0,0){$1$}}%
\put(15.0000,-11.0000){\makebox(0,0){$2$}}%
\put(15.0000,-13.0000){\makebox(0,0){$3$}}%
\put(15.0000,-15.0000){\makebox(0,0){$4$}}%
%
\special{pn 8}%
\special{pa 1600 600}%
\special{pa 1600 1700}%
\special{pa 1600 1700}%
\special{pa 1600 1700}%
\special{fp}%
%
\special{pn 8}%
\special{pa 1400 800}%
\special{pa 2850 800}%
\special{pa 2850 800}%
\special{pa 2850 800}%
\special{fp}%
%
\special{pn 8}%
\special{pa 2800 800}%
\special{pa 2800 1200}%
\special{pa 2350 1200}%
\special{pa 2350 1600}%
\special{pa 1600 1600}%
\special{fp}%
\end{picture}%
\end{center}
Let $J=(x_1^2,x_1x_2,x_1x_3,x_1x_4,x_1x_5,x_2^3,x_2^2x_3,x_2x_3^2,x_3^4)$.
Then the Betti diagram of $J$ is

\begin{center}
\unitlength 0.1in
\begin{picture}( 14.8500, 11.0000)( 13.6500,-17.0000)
\put(17.0000,-9.0000){\makebox(0,0){-}}%
\put(17.0000,-11.0000){\makebox(0,0){$5$}}%
\put(17.0000,-13.0000){\makebox(0,0){$3$}}%
\put(17.0000,-15.0000){\makebox(0,0){$1$}}%
\put(19.5000,-15.0000){\makebox(0,0){$2$}}%
\put(19.5000,-13.0000){\makebox(0,0){$5$}}%
\put(19.5000,-11.0000){\makebox(0,0){$10$}}%
\put(19.5000,-9.0000){\makebox(0,0){-}}%
\put(19.5000,-7.0000){\makebox(0,0){$1$}}%
\put(22.0000,-7.0000){\makebox(0,0){$2$}}%
\put(22.0000,-9.0000){\makebox(0,0){-}}%
\put(22.0000,-11.0000){\makebox(0,0){$10$}}%
\put(22.0000,-13.0000){\makebox(0,0){$2$}}%
\put(22.0000,-15.0000){\makebox(0,0){$1$}}%
\put(24.5000,-15.0000){\makebox(0,0){-}}%
\put(24.5000,-13.0000){\makebox(0,0){-}}%
\put(24.5000,-11.0000){\makebox(0,0){$5$}}%
\put(24.5000,-9.0000){\makebox(0,0){-}}%
\put(24.5000,-7.0000){\makebox(0,0){$3$}}%
\put(27.0000,-7.0000){\makebox(0,0){$4$}}%
\put(27.0000,-9.0000){\makebox(0,0){-}}%
\put(27.0000,-11.0000){\makebox(0,0){$1$}}%
\put(27.0000,-13.0000){\makebox(0,0){-}}%
\put(27.0000,-15.0000){\makebox(0,0){-}}%
\put(17.0000,-7.0000){\makebox(0,0){$0$}}%
\put(15.0000,-9.0000){\makebox(0,0){$1$}}%
\put(15.0000,-11.0000){\makebox(0,0){$2$}}%
\put(15.0000,-13.0000){\makebox(0,0){$3$}}%
\put(15.0000,-15.0000){\makebox(0,0){$4$}}%
%
\special{pn 8}%
\special{pa 1600 600}%
\special{pa 1600 1700}%
\special{fp}%
%
\special{pn 8}%
\special{pa 1400 800}%
\special{pa 2850 800}%
\special{fp}%
%
\special{pn 8}%
\special{pa 2800 800}%
\special{pa 2800 1200}%
\special{pa 2350 1200}%
\special{pa 2350 1600}%
\special{pa 1600 1600}%
\special{fp}%
\end{picture}%
\end{center}
Both $I$ and $J$ admit the extremal area $A$,
and it follows from the Betti diagrams that
$I$ and $J$ have the same Hilbert function.
We claim that there are no graded ideals $L \in \mathcal{L}^A_H$
which satisfies $\beta_{ij}(L) \geq \beta_{ij}(I')$ for all $i,j$
and $I' \in \mathcal{L}^A_H$.

Suppose that there exists such a graded ideal $L \in \mathcal{L}_H^A$.
Then 
$\beta_{2,2+4}(L) \geq \beta_{2,2+4}(I)=2$
and
$\beta_{4,4+2}(L) \geq \beta_{4,4+2}(J)=1$.
On the other hand,
since $L$ have the same Hilbert function as $I$,
we have
$$\sum_{i\geq 0} (-1)^i \beta_{i6}(L) = 
\sum_{i\geq 0} (-1)^i \beta_{i6}(I)=2.$$
However, since $L$ admits $A$,
we have $\sum_{i\geq 0} (-1)^i \beta_{i6}(L)=
\beta_{2,2+4}(L)+ \beta_{4,4+2}(L) \geq 3$.
This is a contradiction.
Thus there exist no graded ideals in $\mathcal{L}_H^A$
which have the maximal graded Betti numbers.
\end{exam}

We will introduce a class of extremal areas $A \subset \area$
which satisfy that, for any numerical function $H:\NN \to \NN$
with $\mathcal{L}^A_H \ne \emptyset$,
there exists a graded ideal $L \in \mathcal{L}^A_H$ such that
$\beta_{ij}(L) \geq  \beta_{ij}(I)$ for all $i,j$ and $I \in\mathcal{L}_H^A$.

\begin{definition}
Let $A \subset \area$ be an extremal area
with the standard representation $A=\bigcup_{k=1}^t \langle (i_k,j_k)\rangle$.
Then $A$ is called a \textit{semi-convex area} if there exists
an integer $1 \leq r \leq t$ such that
$$j_1 =j_2+1 =j_3+2=\cdots =j_r+r-1$$
and
$$i_r=i_{r+1}-1=i_{r+2}-2=\cdots =i_t-(t-r).$$
The element $(i_r,j_r)$ will be called a \textit{top point}
of the semi-convex area $A$.
Note that a top point of a semi-convex area is not always uniquely determined.
Indeed, it is not hard to see that
$(i_r,j_r)$ is a top point of a semi-convex area $A$
if and only if $i_r +j_r= \max \{ i+j: (i,j) \in A \}$.

Let $A \subset \area$ be a semi-convex area.
We say that $(i,j) \in A$ is a \textit{reducible point of $A$}
if $i>0$ and $(i-1,j+1) \not \in A$.
Let $\check A$ denote the subset of $A$
obtained by removing all reducible points from $A$.
\medskip

\begin{center}
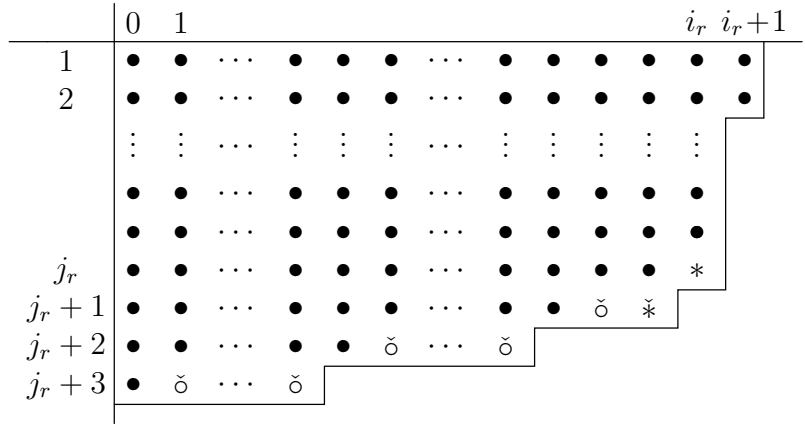
An example of a semi-convex area\medskip\\
\unitlength 0.1in
\begin{picture}( 42.1000, 22.0000)(  5.4000,-26.0000)
\put(12.5000,-5.0000){\makebox(0,0){$0$}}%
\put(15.0000,-7.0000){\makebox(0,0){$\bullet$}}%
\put(15.0000,-9.0000){\makebox(0,0){$\bullet$}}%
\put(15.0000,-11.0000){\makebox(0,0){$\vdots$}}%
\put(15.0000,-14.0000){\makebox(0,0){$\bullet$}}%
\put(15.0000,-16.0000){\makebox(0,0){$\bullet$}}%
\put(15.0000,-18.0000){\makebox(0,0){$\bullet$}}%
\put(15.0000,-20.0000){\makebox(0,0){$\bullet$}}%
\put(15.0000,-22.0000){\makebox(0,0){$\bullet$}}%
\put(15.0000,-24.0000){\makebox(0,0){$\check \circ$}}%
\put(18.0000,-7.0000){\makebox(0,0){$\cdots$}}%
\put(18.0000,-9.0000){\makebox(0,0){$\cdots$}}%
\put(18.0000,-11.5000){\makebox(0,0){$\cdots$}}%
\put(18.0000,-14.0000){\makebox(0,0){$\cdots$}}%
\put(18.0000,-16.0000){\makebox(0,0){$\cdots$}}%
\put(18.0000,-18.0000){\makebox(0,0){$\cdots$}}%
\put(18.0000,-20.0000){\makebox(0,0){$\cdots$}}%
\put(18.0000,-22.0000){\makebox(0,0){$\cdots$}}%
\put(18.0000,-24.0000){\makebox(0,0){$\cdots$}}%
\put(21.0000,-24.0000){\makebox(0,0){$\check \circ$}}%
\put(21.0000,-22.0000){\makebox(0,0){$\bullet$}}%
\put(21.0000,-20.0000){\makebox(0,0){$\bullet$}}%
\put(21.0000,-18.0000){\makebox(0,0){$\bullet$}}%
\put(21.0000,-16.0000){\makebox(0,0){$\bullet$}}%
\put(21.0000,-14.0000){\makebox(0,0){$\bullet$}}%
\put(21.0000,-11.0000){\makebox(0,0){$\vdots$}}%
\put(21.0000,-9.0000){\makebox(0,0){$\bullet$}}%
\put(21.0000,-7.0000){\makebox(0,0){$\bullet$}}%
\put(23.5000,-22.0000){\makebox(0,0){$\bullet$}}%
\put(23.5000,-20.0000){\makebox(0,0){$\bullet$}}%
\put(23.5000,-18.0000){\makebox(0,0){$\bullet$}}%
\put(23.5000,-16.0000){\makebox(0,0){$\bullet$}}%
\put(23.5000,-14.0000){\makebox(0,0){$\bullet$}}%
\put(23.5000,-11.0000){\makebox(0,0){$\vdots$}}%
\put(23.5000,-9.0000){\makebox(0,0){$\bullet$}}%
\put(23.5000,-7.0000){\makebox(0,0){$\bullet$}}%
\put(26.0000,-22.0000){\makebox(0,0){$\check \circ$}}%
\put(26.0000,-20.0000){\makebox(0,0){$\bullet$}}%
\put(26.0000,-18.0000){\makebox(0,0){$\bullet$}}%
\put(26.0000,-16.0000){\makebox(0,0){$\bullet$}}%
\put(26.0000,-14.0000){\makebox(0,0){$\bullet$}}%
\put(26.0000,-11.0000){\makebox(0,0){$\vdots$}}%
\put(26.0000,-9.0000){\makebox(0,0){$\bullet$}}%
\put(26.0000,-7.0000){\makebox(0,0){$\bullet$}}%
\put(29.0000,-7.0000){\makebox(0,0){$\cdots$}}%
\put(29.0000,-9.0000){\makebox(0,0){$\cdots$}}%
\put(29.0000,-11.5000){\makebox(0,0){$\cdots$}}%
\put(29.0000,-14.0000){\makebox(0,0){$\cdots$}}%
\put(29.0000,-16.0000){\makebox(0,0){$\cdots$}}%
\put(29.0000,-18.0000){\makebox(0,0){$\cdots$}}%
\put(29.0000,-20.0000){\makebox(0,0){$\cdots$}}%
\put(29.0000,-22.0000){\makebox(0,0){$\cdots$}}%
\put(32.0000,-22.0000){\makebox(0,0){$\check \circ$}}%
\put(32.0000,-20.0000){\makebox(0,0){$\bullet$}}%
\put(32.0000,-18.0000){\makebox(0,0){$\bullet$}}%
\put(32.0000,-16.0000){\makebox(0,0){$\bullet$}}%
\put(32.0000,-14.0000){\makebox(0,0){$\bullet$}}%
\put(32.0000,-11.0000){\makebox(0,0){$\vdots$}}%
\put(32.0000,-9.0000){\makebox(0,0){$\bullet$}}%
\put(32.0000,-7.0000){\makebox(0,0){$\bullet$}}%
\put(34.5000,-7.0000){\makebox(0,0){$\bullet$}}%
\put(34.5000,-9.0000){\makebox(0,0){$\bullet$}}%
\put(34.5000,-11.0000){\makebox(0,0){$\vdots$}}%
\put(34.5000,-14.0000){\makebox(0,0){$\bullet$}}%
\put(34.5000,-16.0000){\makebox(0,0){$\bullet$}}%
\put(34.5000,-18.0000){\makebox(0,0){$\bullet$}}%
\put(34.5000,-20.0000){\makebox(0,0){$\bullet$}}%
\put(37.0000,-18.0000){\makebox(0,0){$\bullet$}}%
\put(37.0000,-16.0000){\makebox(0,0){$\bullet$}}%
\put(37.0000,-14.0000){\makebox(0,0){$\bullet$}}%
\put(37.0000,-11.0000){\makebox(0,0){$\vdots$}}%
\put(37.0000,-9.0000){\makebox(0,0){$\bullet$}}%
\put(37.0000,-7.0000){\makebox(0,0){$\bullet$}}%
\put(39.5000,-7.0000){\makebox(0,0){$\bullet$}}%
\put(39.5000,-9.0000){\makebox(0,0){$\bullet$}}%
\put(39.5000,-11.0000){\makebox(0,0){$\vdots$}}%
\put(39.5000,-14.0000){\makebox(0,0){$\bullet$}}%
\put(39.5000,-16.0000){\makebox(0,0){$\bullet$}}%
\put(42.0000,-16.0000){\makebox(0,0){$\bullet$}}%
\put(42.0000,-14.0000){\makebox(0,0){$\bullet$}}%
\put(42.0000,-11.0000){\makebox(0,0){$\vdots$}}%
\put(42.0000,-9.0000){\makebox(0,0){$\bullet$}}%
\put(42.0000,-7.0000){\makebox(0,0){$\bullet$}}%
\put(44.5000,-7.0000){\makebox(0,0){$\bullet$}}%
\put(44.5000,-9.0000){\makebox(0,0){$\bullet$}}%
\put(12.5000,-7.0000){\makebox(0,0){$\bullet$}}%
\put(12.5000,-9.0000){\makebox(0,0){$\bullet$}}%
\put(12.5000,-11.0000){\makebox(0,0){$\vdots$}}%
\put(12.5000,-14.0000){\makebox(0,0){$\bullet$}}%
\put(12.5000,-16.0000){\makebox(0,0){$\bullet$}}%
\put(12.5000,-18.0000){\makebox(0,0){$\bullet$}}%
\put(12.5000,-20.0000){\makebox(0,0){$\bullet$}}%
\put(12.5000,-22.0000){\makebox(0,0){$\bullet$}}%
\put(12.5000,-24.0000){\makebox(0,0){$\bullet$}}%
\put(37.0000,-20.0000){\makebox(0,0){$\check \circ$}}%
\put(39.5000,-20.0000){\makebox(0,0){$\check *$}}%
\put(39.5000,-18.0000){\makebox(0,0){$\bullet$}}%
\put(42.0000,-18.0000){\makebox(0,0){$*$}}%
\put(42.0000,-16.0000){\makebox(0,0){$\bullet$}}%
\put(15.0000,-5.0000){\makebox(0,0){$1$}}%
\put(45.0000,-5.0000){\makebox(0,0){$i_r\! + \! 1$}}%
\put(42.0000,-5.0000){\makebox(0,0){$i_r$}}%
\put(9.0000,-7.0000){\makebox(0,0){$1$}}%
\put(9.0000,-9.0000){\makebox(0,0){$2$}}%
\put(9.0000,-24.0000){\makebox(0,0){$j_r +3$}}%
\put(9.0000,-22.0000){\makebox(0,0){$j_r +2$}}%
\put(9.0000,-20.0000){\makebox(0,0){$j_r +1$}}%
\put(9.0000,-18.0000){\makebox(0,0){$j_r$}}%
%
\special{pn 8}%
\special{pa 800 600}%
\special{pa 4750 600}%
\special{pa 4750 600}%
\special{pa 4750 600}%
\special{fp}%
%
\special{pn 8}%
\special{pa 1150 400}%
\special{pa 1150 2600}%
\special{fp}%
%
\special{pn 8}%
\special{pa 600 600}%
\special{pa 800 600}%
\special{fp}%
%
\special{pn 8}%
\special{pa 4550 600}%
\special{pa 4550 1000}%
\special{pa 4350 1000}%
\special{pa 4350 1900}%
\special{pa 4100 1900}%
\special{pa 4100 2100}%
\special{pa 3350 2100}%
\special{pa 3350 2300}%
\special{pa 2250 2300}%
\special{pa 2250 2500}%
\special{pa 1150 2500}%
\special{fp}%
\end{picture}%
\end{center}
(In the picture, $*$ and $\check *$ are top points.
Also, $\check \circ$ and $\check *$ are reducible points.)
\end{definition}

We will show that if $ A \subset \area$ is a semi-convex area
then, for any numerical function $H:\NN \to \NN$ with $\mathcal{L}_H^A\ne \emptyset$,
there exists a monomial ideal $L \in \mathcal{L}_H^A$ which has the maximal graded Betti numbers among graded ideals in $\mathcal{L}_H^A$.

Let $I \subset S$ be a strongly stable ideal.
For integers $1 \leq i \leq n$ and $j \geq 0$,
define
$$M_i(I,j)=\{u \in \M_j: u \in I \mbox{ and } \max(u)=i\}$$
and 
$$M_{\leq i}(I,j)=\{u \in \M_j: u \in I \mbox{ and } \max(u) \leq i\},$$
where $\max(1)=0$.

\begin{lemma}\label{shadow}
Let $I \subset S$ be a strongly stable ideal.
Then $|M_{i}(I,j)| \geq |M_{\leq i}(I,j-1)|$ for all $ 1 \leq i \leq n$
and $j > 0$.
Moreover, if $\beta_{ii+j}(I)=0$
then $|M_{i+1}(I,j)| = |M_{\leq i+1}(I,j-1)|$.
\end{lemma}

\begin{proof}
Since $x_{i}M_{\leq i}(I,j-1) \subset M_{i}(I,j)$,
the first statement is obvious.
On the other hand, if $\beta_{ii+j}(I)=0$ then the Eliahou--Kervaire
formula says that $I$ has no generators $u \in G(I)$ of degree $j$
with $\max(u)\geq i+1$.
Then, since $I$ is strongly stable,
it follows from \cite[Lemma 2.9]{H} that 
\begin{eqnarray*}
\ M_{i+1}(I,j) &=&\{x_ku:k=1,\dots,i+1,\ u \in M_{\leq i+1}(I,j-1),\ \max(ux_k)=i+1\}\\
&=& x_{i+1}  M_{\leq i+1}(I,j-1).
\end{eqnarray*}
Thus the claim follows.
\end{proof}

\begin{lemma} \label{carea}
Let $A \subset \area$ be a semi-convex area with a top point $(i_r,j_r) \in A$.
\begin{itemize}
\item[(i)] If $(i,j) \not \in A$  and $j \leq j_r$ then $(i+1,j-1) \not \in A$;
\item[(ii)] If $(i,j) \not \in A$  and $j \geq j_r$
then $(i-1,j+1) \not \in A$;
\item[(iii)] If $(i,j) \in A$ is reducible then $j \geq j_r$.
\end{itemize}
\end{lemma}

\begin{proof}
Let $A =\bigcup_{k=1}^t \langle (i_k,j_k) \rangle$ be the standard representation of $A$.

(i) It is enough to show that if $(i+1,j-1) \in A$ and $j \leq j_r$ then $(i,j) \in A$.
Suppose $(i+1,j-1) \in A$ and $j \leq j_r$.
If $(i+1,j-1) \in \langle (i_r,j_r) \rangle$ then $i<i_r$.
Since $j \leq j_r$ we have
 $(i,j) \in \langle (i_r,j_r) \rangle \subset A$.
Assume that $(i+1,j-1)  \in A \setminus\langle (i_r,j_r) \rangle$.
Then, since $j \leq j_r$, there exists an integer $k>r$ such that $(i+1,j-1) \in \langle (i_k,j_k) \rangle$.
The definition of semi-convex areas says that $i_{k-1}+1 = i_k$.
Also, we have $j_{k-1}>j_k$ by the definition of standard representations.
Then, $i+1 \leq i_k =i_{k-1} + 1$ and $j-1 \leq j_k \leq j_{k-1}-1$.
Thus we have $(i,j) \in \langle (i_{k-1},j_{k-1}) \rangle \subset A$.

(ii)
It suffices to show that if $(i-1,j+1)\in A$ and $j \geq j_r$ then
$(i,j)\in A$.
Suppose $(i-1,j+1) \in A$ and $j \geq j_r$.
Since $j \geq j_r$, there exists a $k<r$ such that
$(i-1,j+1) \in \langle (i_k,j_k) \rangle$.
The definition of semi-convex areas says that
$j+1 \leq j_k =j_{k+1}+1$.
Also, we have
$i \leq i_k +1 \leq i_{k+1}$ by the definition of standard representations.
Thus we have $(i,j) \in \langle (i_{k+1},j_{k+1})\rangle \subset A$ as desired.

(iii)
If $(i,j) \in A$ is a reducible point of $A$ then $i>0$ and $(i-1,j+1) \not \in A$.
However, statement (i) says that $(i-1,j+1) \not \in A$ and $j < j_r$
imply $(i,j) \not \in A$.
Hence $j \geq j_r$.
\end{proof}

\begin{lemma} \label{choice}
Let $A \subset \area$ be a semi-convex area,
and let $I \subset S$ and $J \subset S$ be strongly stable ideals
which admit $A$.
If $I$ and $J$ have the same Hilbert function,
then we have
\begin{eqnarray}
|M_{i+1}(I,j)|=|M_{i+1}(J,j)| \ \ \mbox{ for all } (i,j) \not \in \check A. \label{CC1}
\end{eqnarray}
\end{lemma}

\begin{proof}
Let $A =\bigcup_{k=1}^t \langle (i_k,j_k) \rangle$ be the standard representation of $A$
and $(i_r,j_r) \in A$ a top point of $A$.
$|M_1(I,j)|=|M_1(J,j)|$ for all $j \geq 0$ is obvious.
Thus we assume $i>0$.
\medskip

\textbf{[Case 1]}
First, we consider the case $j<j_r$.
We use induction on $j$.
Since $I_1 =J_1$, the statement is obvious for $j=1$.
We assume that $j<j_r$ and (\ref{CC1}) holds for all $(i,k) \not \in \check A$ with $k<j$.

Let $(i,j) \not \in \check A$.
Since $j<j_r$, Lemma \ref{carea} (iii) says that
$(i,j)$ is not a reducible point.
Hence $(i,j) \not \in A$.
Since $I$ and $J$ admit $A$,
we have $\beta_{ii+j}(I)=\beta_{ii+j}(J)=0$.
Then Lemma \ref{shadow} says that
\begin{eqnarray}
|M_{i+1}(I,j)|=|M_{\leq i+1}(I,j-1)|=\dim_K I_{j-1} - \sum_{k=i+2}^{n}|M_k(I,j-1)|
\label{ML1}
\end{eqnarray}
and
\begin{eqnarray}
|M_{i+1}(J,j)|=|M_{\leq i+1}(J,j-1)|=\dim_K J_{j-1} - \sum_{k=i+2}^{n} |M_k(J,j-1)|.
\label{ML2}
\end{eqnarray}
By Lemma \ref{carea} (i), we have $(i+1,j-1) \not \in A$.
Then, since $A$ is an extremal area, we have $(p,j-1) \not \in A$ for all $p \geq i+1$.
Thus we have $|M_{i+1}(I,j)|=|M_{i+1}(J,j)|$ by
(\ref{ML1}) and (\ref{ML2}) together with the induction hypothesis.
\medskip

\textbf{[Case 2]}
Second, we consider the case $j \geq j_r$.

If $j  \geq j_1$,
then Lemma \ref{ek} says that $I_{\geq j}$ and $J_{\geq j}$ are
strongly stable ideals generated in degree $j$
since $I$ and $J$ admit $A$.
Since $I_{\geq j}$ and $J_{\geq j}$ have the same Hilbert function,
Lemma \ref{equive} says
$$|M_i(I,j)|=\ell_i(I_{\geq j})=\ell_i(J_{\geq j})=|M_i(I,j)|
 \ \ \mbox{ for all } i.$$

Next, we will show the statement for $j_r \leq j < j_1$
by using induction on $j$.
Assume that $j_r \leq j <j_1$ and (\ref{CC1}) holds for all $(i,k) \not \in \check A$ with $k>j$.

Let $(i,j) \not \in \check A$.
Then we have $(i',j) \not \in \check A$ for all $i' \geq i$,
and Lemma \ref{carea} (ii) together with the definition of reducible points says that $(i'-1,j+1) \not \in A$
for all $i' \geq i$.
Thus $\beta_{i'-1,i'+j}(I)=\beta_{i'-1,i'+j}(J)=0$ for all
$i' \geq i$.
Then Lemma \ref{shadow} says that, for all $i' \geq i$, we have
\begin{eqnarray}
|M_{i'}(I,j+1)|=|M_{\leq i'}(I,j)|= \dim_K I_j - \sum_{k=i'+1}^n
|M_k(I,j)|
\label{ML3}
\end{eqnarray}
and
\begin{eqnarray}
|M_{i'}(J,j+1)|=|M_{\leq i'}(J,j)|= \dim_K J_j - \sum_{k=i'+1}^n
|M_k(J,j)|.
\label{ML4}
\end{eqnarray}
Notice that the induction hypothesis says $|M_{i'}(I,j+1)|=|M_{i'}(J,j+1)|$ for $i' \geq i$ since $(i'-1,j+1) \not \in A$.
Then $|M_n(I,j)|=|M_n(J,j)|$ immediately follows from the above equations
together with the induction hypothesis.
By arguing inductively,
(\ref{ML3}) and (\ref{ML4}) together with the induction hypothesis
imply $|M_{i'+1}(I,j)|=|M_{i'+1}(J,j)|$ for all $i' \geq i$.
In particular, we have $|M_{i+1}(I,j)|=|M_{i+1}(J,j)|$ for all $i$
with $(i,j) \not \in \check A$ as desired.
\end{proof}

To construct an ideal which gives the maximal graded Betti
numbers in $\mathcal{L}_H^A$,
we need to introduce the following lemma.

\begin{lemma} \label{mvector}
Let $d$ be a positive integer.
If $V \subset \M_d$ is strongly stable then,
for any integer $1 < r \leq n+1$,
there exists the unique $d$-linear lexsegment subset $W \subset \M_d$ which satisfies the following conditions.
\begin{itemize}
\item[(i)] $\ell_k(V)=\ell_k(W)$ for all $k \geq r$;
\item[(ii)] $|M_{\leq r-1}(W)|=|M_{\leq r-1}(V)|$ and
$M_{\leq r-1}(W)$ is lexsegment in $K[x_1,\dots,x_{r-1}]$.
\end{itemize} 
\end{lemma}

\begin{proof}
Let $W_1 \subset \M_d \cap K[x_1,\dots,x_{r-1}]$ be the lexsegment set of monomials with $|W_1|=|M_{\leq r-1}(V)|$.
The vectors $\ell(W_1)=(\ell_1(W_1),\dots,\ell_{r-1}(W_1))$
and $\ell(V)=(\ell_1(V),\dots,\ell_n(V))$
are $M$-vectors satisfying
$\ell_2 (W_1) \leq d$ and $\ell_2 (V) \leq d$ by Proposition \ref{lsequence}.
On the other hand,
Lemma \ref{bigatti} says that
\begin{eqnarray*}
\ell_{r-1}(W_1)=|M_{ r-1}(W_1)|&=&|W_1|-|M_{\leq r-2}(W_1)|\\
&\geq& |M_{\leq r-1}(V)| -|M_{\leq r-2}(V)|\\
& =& |M_{ r-1}(V)|=\ell_{r-1}(V).
\end{eqnarray*}
Then $\ell=(\ell_1(W_1),\dots,\ell_{r-1}(W_1),
\ell_r(V),\dots,\ell_n(V))$ is also an $M$-vector.
Thus, by Propositions \ref{llex} and \ref{lsequence},
there exists the $d$-linear lexsegment subset $W \subset \M_d$
with $\ell(W) = \ell$.
This set $W$ satisfies condition (i).

On the other hand, for $k<r$,
both $D_k(W)$ and 
$D_k(W_1)$ are lexsegment sets of monomials in $K[x_1,\dots,x_k]$.
Since $|D_k(W)|=\ell_k(W)=\ell_k(W_1)=|D_k(W_1)|$,
we have $D_k(W)=D_k(W_1)$ for all $k <r$.
Thus $\M_{\leq r-1}(W)=\bigcup_{k=1}^{r-1} x_kD_k(W)= W_1$ is lexsegment in $K[x_1,\dots,x_{r-1}]$.
Hence $W$ satisfies condition (ii).

Recall that $d$-linear lexsegment subsets are uniquely determined from
their $\ell$-sequence.
Since $\ell$ is uniquely determined from $V$ and $r$,
this set $W$ is uniquely determined from $V$ and $r$.
\end{proof}

\begin{construction} \label{defi}
Let $A \subset \area$ be a semi-convex area,
$A= \bigcup_{k=1}^t \langle (i_k,j_k) \rangle$ the standard representation of $A$,
$(i_r,j_r) \in A$ a top point of $A$
and $H:\NN \to \NN$ a numerical function with $\mathcal{L}_H^A \ne \emptyset$.
Set $p_j=\max\{i: (i,j) \in A\}$ for $j=1,2,\dots,j_1$ and $p_j=-1$ for $j >j_1$.

We will construct an ideal which gives the maximal graded Betti numbers
among graded ideals in $\mathcal{L}_H^A$ from a strongly stable ideal in
$\mathcal{L}_H^A$.
Since we assume $\mathrm{char}(K)=0$,
Lemma \ref{bcp} says that
if $\mathcal{L}_H^A \ne \emptyset$ then there exists a strongly stable ideal
$I \in \mathcal{L}_H^A$.
We construct the monomial ideal $\Lex(I,A) \subset S$ from $I$ as follows:
\medskip

(A)
For $j<j_r$,
let $L_j(I,A) \subset \M_j \cap K[x_1,\dots,x_{p_j+1}]$
be the lexsegment set of monomials in $K[x_1,\dots,x_{p_j+1}]$
with $|L_j(I,A)|=|M_{\leq p_j+1}(I,j)|$.

(B)
For $j_r \leq j \leq j_1$,
since $M_{\leq p_j+1}(I,j) \subset K[x_1,\dots,x_{p_j+1}]$
is a strongly stable set of monomials and $p_{j+1}+1 \leq p_j$
by the definition of semi-convex areas,
Lemma \ref{mvector} says that there exists the $j$-linear lexsegment subset $L_j \subset \M_j \cap K[x_1,\dots,x_{p_j+1}]$
satisfying
\begin{itemize}
\item[(i)]
$|M_{i+1}(L_j)|=|M_{i+1}(I,j)|$ for all $p_j \geq i > p_{j+1}+1$;
\item[(ii)] $M_{\leq p_{j+1}+2}(L_j)$ is the lexsegment set of monomial in $K[x_1,\dots,x_{p_{j+1}+2}]$ with 
$|M_{\leq p_{j+1}+2}(L_j)| = |M_{\leq p_{j+1}+2}(I,j)|$.
\end{itemize}
Set $L_j(I,A)=L_j$.
\medskip

Define
$$\Lex(I,A)= \sum_{j=1}^{j_1} \langle L_j(I,A) \rangle, $$
where $\langle L_j(I,A) \rangle$ is the ideal generated by the set of monomials $L_j(I,A)$.
\end{construction}

We show that $\Lex(I,A)$ admits $A$,
and only depends on the semi-convex area $A$ and the Hilbert function
of $I$.

\begin{lemma} \label{independent}
With the same notations as in Construction \ref{defi},
\begin{itemize}
\item[(i)] $\Lex(I,A)$ is independent of the choice of a top point $(i_r,j_r) \in A$ of
$A$;
\item[(ii)] $\Lex(I,A)$ is independent of the choice of a strongly stable ideal $I \in \mathcal{L}_H^A$.
\end{itemize}
\end{lemma}

\begin{proof}
(i)
Assume that $A$ has two top points $(i_r,j_r) \in A$ and $(i_{r'},j_{r'}) \in A$ with $r<r'$.
Then since both elements are top points of $A$,
we have
$i_r=i_{r+1}-1 =i_{r+2}-2 =\cdots = i_{r'}-(r'-r)$
and 
$j_r=j_{r+1}+1 =j_{r+2}+2 =\cdots = j_{r'}+(r'-r)$.

This property says that $p_j =p_{j+1}+1$ for $j_{r'} \leq j < j_{r}$.
Then construction (A) is the same as construction (B) for $j_{r'} \leq j < j_{r}$.
Thus the claim follows.

(ii) Let $(i_r,j_r) \in A$ be a top point of $A$.
For $j<j_r$, construction (A) says that $L_j(I,A)$ only depends on the number
$|M_{\leq p_j+1}(I,j)| = \dim_K I_j - \sum_{i =p_j+2}^{n} |M_{i}(I,j)|$.
However, if $i>p_j$ then $(i,j) \not \in A$.
Thus Lemma \ref{choice} says that
this number is independent of the choice of a strongly stable ideal $I \in \mathcal{L}_H^A$.

For $j_1 \geq j \geq j_r$,
the set $L_j(I,A)$ only depends on $|M_{\leq p_{j+1}+2}(I,j)|$ together with 
$|M_{i+1}(I,j)|$ for $p_j \geq i>p_{j+1} +1$.
However, if $i>p_{j+1}+1$ then $(i,j) \not \in \check A$.
Thus Lemma \ref{choice} says that 
these numbers are also independent of the choice of a strongly stable ideal $I \in \mathcal{L}_H^A$.
\end{proof}

\begin{lemma} \label{admit}
With the same notations as in Construction \ref{defi},
$\Lex(I,A)$ is a strongly stable ideal which admits $A$.
\end{lemma}

\begin{proof}
A sum of strongly stable ideals is again strongly stable.
Thus $\Lex(I,A)$ is strongly stable.
On the other hand, since $L_j(I,A) \subset K[x_1,\dots,x_{p_j+1}]$,
it follows that $\Lex(I,A)$ has no generators $u \in G(\Lex(I,A))$ of degree $j$ with $\max(u) >p_j+1$.
Then the Eliahou--Kervaire formula says $\beta_{ii+j}(\Lex(I,A))=0$
for all $(i,j) \not \in A$.
Hence $\Lex(I,A)$ admits $A$.
\end{proof}

The next lemma immediately follows from the definition of $\Lex(I,A)$ and Lemma \ref{bigatti}.

\begin{lemma} \label{convbig}
With the same notations as in Construction \ref{defi},
for each $j=1,\dots,j_1$,
$L_j(I,A) \subset K[x_1,\dots,x_{p_j+1}]$ is a $j$-linear lexsegment set of monomials
satisfying
\begin{itemize}
\item[(i)] $|L_j(I,A)|=|M_{\leq p_j+1} (I,j)|$;
\item[(ii)] $|M_{\leq i+1}(L_j(I,A))| \leq |M_{\leq i+1}(I,j)|$ for all $i \leq p_j$.
\end{itemize}
\end{lemma}

Next, we will show $I$ and $\Lex(I,A)$ have the same Hilbert function.

\begin{lemma} \label{shade}
With the same notations as in Construction \ref{defi},
let $W_j$ be the set of monomials in $\Lex(I,A)$ of degree $j$
and $V_j$ the set of monomials in $I$ of degree $j$
for all $j \geq 0$.
Then
\begin{itemize}
\item[(a)] $M_{\leq p_j +1}(W_j)=L_j(I,A)$ for all $j \leq j_1$;
\item[(b)] $|W_j|=|V_j|$ for all $j \geq 0$;
\item[(c)] $|M_{i+1}(W_j)|=|M_{i+1}(V_j)|$ for all $(i,j) \not \in \check A$;
\item[(d)] $|M_{\leq i}(W_j)|=|M_{\leq i}(V_j)|$ for all $(i,j) \not \in \check A$.
\end{itemize}
\end{lemma}

\begin{proof}
To simplify the argument, let $L_j=L_j(I,A)$ for all $j \leq j_1$.
We use induction on $j$.
For $j=1$, since $W_1 =L_1 =V_1$, the assertion is obvious.

Let $d >1$ be an integer.
Assume that the statements are true for $j <d$.
We will prove the statements for $j=d$.

First, we claim that
\begin{eqnarray}
&&|M_{i+1}(V_d)|= |M_{\leq i+1}(V_{d-1})| =
|M_{\leq i+1}(W_{d-1})|=|M_{i+1}(W_{d})| \ \mbox{ for all } i>p_d.
\label{LLA}
\end{eqnarray}
Since $I$ and $\Lex(I,A)$ admit $A$,
we have $\beta_{ii+d}(I)= \beta_{ii+d}(\Lex(I,A))=0$
for all $i > p_d$.
Thus the first equality of (\ref{LLA}) and the third one
follow from Lemma \ref{shadow}.
On the other hand, since $i>p_d$ implies $(i,d) \not \in A$ and
$(i,d) \not \in A$ implies $(i+1,d-1) \not \in \check A$ by the definition of reducible points,
the second equality of (\ref{LLA}) follows from the induction hypothesis.

Second, we will show statement (a).
Clearly, we have
$$M_{\leq p_d +1} (W_d)= L_d \cup \{x_k u : u\in M_{\leq p_d +1}(W_{d-1})
,\ k=1,2,\dots,p_d+1\}.$$
Since $p_{d-1} \geq p_d$, the induction hypothesis says that
$M_{\leq p_d +1}(W_{d-1})= M_{\leq p_d +1}(L_{d-1})$.
Then since $L_d$ is strongly stable,
to prove statement (a), it is enough to prove that
\begin{eqnarray}
x_{p_d +1} M_{\leq p_d +1}(L_{d-1}) \subset L_d.
\label{LLW}
\end{eqnarray}
By Lemmas \ref{shadow} and \ref{convbig}, we have
$$ |M_{\leq p_d +1}(L_{d-1})| \leq |M_{\leq p_d +1}(V_{d-1})|
\leq |M_{p_d +1}(V_d)|.$$
Lemma \ref{convbig} also says that
\begin{eqnarray*}
|M_{p_d +1}(V_d)| &=& |M_{\leq p_d+1}(V_d)| - |M_{\leq p_d}(V_d)|\\
&\leq& |L_d| - |M_{\leq p_d}(L_d)|
= |M_{p_d+1}(L_d)|
= |D_{p_d+1}(L_d)|,
\end{eqnarray*}
where $D_{p_d+1}(L_d)=\{ (u/x_{p_d +1}): u \in L_d,\ \max(u) =p_d +1\}$.
Thus $|M_{\leq p_d +1}(L_{d-1})| \leq |D_{p_d+1}(L_d)|$.
Since $L_d$ is $d$-linear lexsegment,
$D_{p_d +1}(L_d) \subset \M_{d-1} \cap K[x_1,\dots,x_{p_d +1}]$
is lexsegment in $K[x_1,\dots,x_{p_d+1}]$.
On the other hand, the construction of $L_{d-1}$ says that $M_{\leq p_d +1}(L_{d-1}) \subset \M_{d-1} \cap K[x_1,\dots,x_{p_d +1}]$
is also lexsegment in $K[x_1,\dots,x_{p_d+1}]$.
Thus we have 
$$x_{p_d +1} M_{\leq p_d +1}(L_{d-1}) \subset
x_{p_d +1} D_{p_d+1}(L_d) \subset L_d.$$
This is (\ref{LLW}).
Hence statement (a) follows.

Since we proved $M_{\leq p_d +1}(W_d)= L_d$,
Lemma \ref{convbig} (i) says  $|M_{\leq p_d +1}(W_d)|= |L_d|=|M_{\leq p_d +1}(V_d)|$.
Then, by using (\ref{LLA}),
we have
\begin{eqnarray*}
\ |W_d| = |M_{\leq p_d +1 }(W_d)| + \sum_{i=p_d +2}^n |M_i(W_d)|
= |M_{\leq p_d +1 }(V_d)| + \sum_{i=p_d +2}^n |M_i(V_d)| = |V_d|.
\end{eqnarray*}
Thus statement (b) follows.

Next, we will prove statement (c).
If $(i,d) \not \in A$ then $i>p_d$.
Thus statement (c) is equal to (\ref{LLA}) in this case.
Also, if $i=0$ then statement (c) is obvious.
Suppose that $(i,d) \in A$ but $(i,d) \not \in \check A$ and $i>0$.
Then $(i,d)$ is a reducible point of $A$.
Thus $(i-1,d+1) \not \in A$ and $(i,d) \in A$.
Hence $p_{d+1} +1 < i \leq p_d$.
Let $(i_r,j_r) \in A$ be a top point of $A$.
Then, we have $d \geq j_r$ by Lemma \ref{carea} (iii).
Then construction (B) of $L_d(I,A)$ says
$|M_{i+1}(L_d)|=|M_{i+1}(V_d)|$.
Since we already proved $L_d =M_{\leq p_d +1}(W_d)$,
we have
$|M_{i+1}(W_d)| =|M_{i+1}(L_d)| =|M_{i+1}(V_d)|$ as desired.
Thus statement (c) follows.

Finally,
it remains to show statement (d).
However, statement (b) and (c) imply that,
for all $(i,d) \not \in \check A$, we have
\begin{eqnarray*}
|M_{\leq i}(W_d)| = |W_d| - \sum_{k=i+1}^n |M_{k}(W_d)|
= |V_d| - \sum_{k=i+1}^n |M_{k}(V_d)| = |M_{\leq i}(V_d)|,
\end{eqnarray*}
as required.
\end{proof}

\begin{proposition} \label{main7}
With the same notation as in Construction \ref{defi},
\begin{itemize}
\item[(i)] $I$ and $\Lex(I,A)$ have the same Hilbert function;
\item[(ii)] $|M_{\leq i}(I,j)| \geq |M_{\leq i}(\Lex(I,A),j)|$
for all $i$ and $j$; 
\item[(iii)] $\beta_{ij}(I) \leq \beta_{ij}(\Lex(I,A))$ for all $i$ and $j$.
\end{itemize}
\end{proposition}

\begin{proof}
The first statement is just Lemma \ref{shade} (b).
The second statement follows from Lemma \ref{convbig} and
Lemma \ref{shade} (a) and (d).
The last statement follows from the second statement by using
(\ref{LRM6}).
\end{proof}

\begin{theorem} \label{laststand}
Assume that $\mathrm{char}(K)=0$.
Let $A \subset \area$ be a semi-convex area and
$H : \NN \to \NN$ a numerical function with $\mathcal{L}_H^A \ne \emptyset$.
Then there exists a strongly stable ideal $L \in \mathcal{L}_H^A$ such that
$\beta_{ij}(L) \geq \beta_{ij}(I)$ for all $i,j$ and $I \in \mathcal{L}_H^A$.
\end{theorem}

\begin{proof}
By Lemmas \ref{itumono} and \ref{bcp}, for any graded ideal $I \in \mathcal{L}_H^A$,
there exists a strongly stable ideal $J \in \mathcal{L}_H^A$ satisfying $\beta_{ij}(I) \leq \beta_{ij}(J)$ for all $i$ and $j$.
On the other hand,
Proposition \ref{main7} says that
$\beta_{ij}(J) \leq \beta_{ij}(\Lex(J,A))$ for all $i$ and $j$,
and Lemma \ref{independent} says that the ideal $\Lex(J,A)$ is uniquely determined from $A$ and $H$.
Since Lemma \ref{admit} and Proposition \ref{main7}
say $\Lex(J,A) \in \mathcal{L}_H^A$,
this ideal $\Lex(J,A)$ satisfies the required condition.
\end{proof}

Let $H:\NN \to \NN$ be a numerical function.
For integers $1 \leq p \leq n$ and $d>0$,
let $\mathcal{L}^{(p,d)}_H$ be the set of graded ideals defined by
$$ \mathcal{L}^{(p,d)}_H=
\{ I \subset S: \proj (S/I) \leq p,\ \reg(I) \leq d,\ H(S/I,t)=H(t) \mbox{ for all } t \in \NN\}.$$
It is clear that if $A=\langle (i,j) \rangle$ for some $(i,j) \in \area$
then $A$ is a semi-convex area.
Since $\mathcal{L}^{(p,d)}_H=\mathcal{L}^{\langle (p-1,d) \rangle}_H$,
the following result immediately follows from Theorem \ref{laststand}.

\begin{cor}
Assume that $\mathrm{char}(K)=0$.
Let  $1 \leq p \leq n$ and $d>0$ be integers.
If $H:\NN \to \NN $ is a numerical function with $\mathcal{L}^{(p,d)}_H \ne \emptyset$
then there exists a strongly stable ideal $L \in \mathcal{L}^{(p,d)}_H$
such that $\beta_{ij}(L) \geq \beta_{ij}(I)$ for all $i,j$ and $I \in \mathcal{L}^{(p,d)}_H$. 
\end{cor}

Next, we consider squarefree monomial ideals.
Let $A \subset \area$ be an extremal area
and $H:\NN \to \NN$ a numerical function.
We write
$$\mathcal{SL}_H^A =\{ I \subset S: I \in \mathrm{Sq}(S),\ I \mbox{ admits }A \mbox{ and }
H(S/I,t)=H(t) \mbox{ for all } t \in \NN\},$$
where $\mathrm{Sq}(S)$ is the set of all squarefree monomials in $S$.
Set
$Q=\{(i,j) \in \area: i+j \leq n\}$.
Then it is easy to see that if $A$ is a semi-convex area
then $A\cap Q$ is again a semi-convex area.
Also since any squarefree monomial ideal admits $Q$,
we have $\mathcal{SL}_H^A = \mathcal{SL}_H^{A \cap Q}$.

\begin{cor} \label{specialsq}
Assume that $\mathrm{char}(K)=0$.
Let $A \subset \area$ be a semi-convex area and
$H:\NN \to \NN$ a numerical function with $\mathcal{SL}_H^A \ne \emptyset$.
Then there exists a squarefree strongly stable ideal $L \in \mathcal{SL}_H^A$
such that $\beta_{ij}(L) \geq \beta_{ij}(I)$ for all $i,j$ and $I \in \mathcal{SL}_H^A$. 
\end{cor}

\begin{proof}
Since $A \cap Q$ is a semi-convex area, there exists a strongly stable ideal $J \in \mathcal{L}_H^{A \cap Q}$
such that $\beta_{ij}(J) \geq \beta_{ij}(I)$ for all $i,j$ and $I \in \mathcal{L}_H^{A \cap Q}$.
Since $J$ admit $Q$, $J$ satisfies the assumption of Lemma \ref{selfoperator}.
Set $L= \tilde \Phi(J)$.
Then Lemma \ref{selfoperator} says that $L \in \mathcal{SL}_H^{A}$ and
$\beta_{ij}(L)= \beta_{ij}(J)$ for all $i$ and $j$.
Since $\mathcal{L}_H^{A \cap Q} \supset \mathcal{SL}_H^{A \cap Q} = \mathcal{SL}_H^{A}$,
this ideal $L$ satisfies the required condition.
\end{proof}

Let $H:\NN \to \NN$ be a numerical function.
For integers $1 \leq p \leq n$ and $0 \leq d \leq n$,
let 
$ \mathcal{SL}^{(p,d)}_H= \mathcal{L}^{(p,d)}_H \cap \mathrm{Sq}(S),
$
in other words, $\mathcal{SL}^{(p,d)}_H= \mathcal{SL}^{\langle (p-1,d) \rangle}_H$.
Then Corollary \ref{specialsq} implies

\begin{cor}
\label{wasure}
Assume that $\mathrm{char}(K)=0$.
Let  $1 \leq p \leq n$ and $1 \leq d \leq n$ be integers.
If $H:\NN \to \NN $ be a numerical function with $\mathcal{SL}^{(p,d)}_H \ne \emptyset$
then there exists a squarefree strongly stable ideal $L \in \mathcal{SL}^{(p,d)}_H$
such that $\beta_{ij}(L) \geq \beta_{ij}(I)$ for all $i,j$ and $I \in \mathcal{SL}^{(p,d)}_H$. 
\end{cor}

\begin{exam}
Let $S=K[x_1,\dots,x_5]$,
$I=(x_1^2,x_1x_2,x_1x_3,x_1x_4,x_2^2,x_2x_3^3,x_3^4)$
and $H$ the Hilbert function of $S/I$.
Set $A=\langle (2,4) \rangle \cup \langle (4,2) \rangle$.
In Example \ref{counter},
it was proved that there are no graded ideals having the maximal graded Betti numbers among all graded ideals in $\mathcal{L}_H^A$.
However, if we consider the semi-convex area
$B=\langle (2,4) \rangle\cup \langle (3,3) \rangle \cup \langle (4,2) \rangle \supset A,$
$\Lex(I,B)$ has the maximal graded Betti numbers among graded ideals in $\mathcal{L}_H^B$.
Indeed,
$$\Lex(I,B) = (x_1^2,x_1x_2,x_1x_3,x_1x_4,x_1x_5,x_2^3,x_2^2x_3,x_2^2x_4,
x_2x_3^3,x_3^4)$$
and its Betti diagram is
\begin{center}
\begin{tabular}{cc|ccccccccc}
& & 0 & 1 & 2 & 3 & 4\\
\hline
&2 & 5 & 10 & 10 & 5 & 1 &\\
&3 & 3 & 6 & 4 & 1 & - \\
&4 & 2 & 4 & 2 &- & -\\
&&\vspace{-5pt}\\
& total & 10 & 20 & 16 & 6 & 1
\end{tabular}\medskip
\end{center}
Remark that the graded Betti numbers of $\Lex(I)$ are much larger than those of $\Lex(I,B)$.
Indeed, we have $\reg(\Lex(I))=17$ and $|G(\Lex(I))|=38$.
\end{exam}

The above example leads us to ask
whether there exists the smallest semi-convex area containing
an extremal area $A$.
Constructing such semi-convex areas is not difficult.
We can construct them as follows:

Let $A \subset \area$ be an extremal area with the standard representation
$A=\bigcup_{k=1}^t \langle (i_k,j_k)\rangle$.
Choose any point $(i_r,j_r)$ satisfying $i_r + j_r = \max\{i+j:(i,j) \in A\}$.
Let $\mathrm{conv}(A)$ be the extremal area defined by
\begin{eqnarray*}
\mathrm{conv}(A)=&&
\bigg(\bigcup_{k=1}^{r-1} \bigg\{ \bigcup_{p=0}^{j_k-j_r-1}\langle (i_k +p,j_k-p)\rangle\bigg\}\bigg)
 \\
&&
 \bigcup \bigg(\bigcup_{k=r+1}^t \bigg\{ \bigcup_{p=0}^{i_k-i_r-1}\langle (i_k -p,j_k+p)\rangle\bigg\} \bigg) \cup \langle (i_r,j_r)\rangle.
\end{eqnarray*}
Then it is not hard to see that $\mathrm{conv(A)}$ is independent of the choice of a point
$(i_r,j_r) \in A$ with $i_r +j_r =\max\{i+j : (i,j) \in A\}$,
and is indeed the unique smallest semi-convex area containing $A$.
The semi-convex area $\mathrm{conv}(A) \subset \area $ 
will be called the \textit{semi-convex hull of $A$.}

\begin{center}
The semi-convex hull of $A= \cup_{k-1}^4 \langle
(i_k,j_k)\rangle$ \medskip\\
\unitlength 0.1in
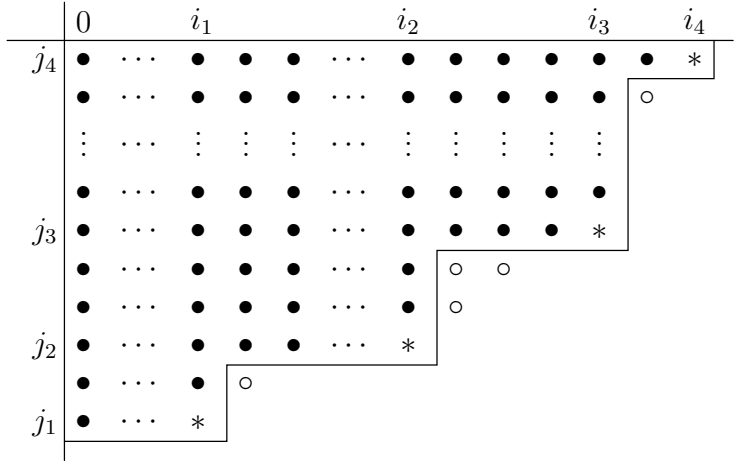
\begin{picture}( 38.2500, 24.0000)( 10.7500,-28.0000)
\put(15.0000,-5.0000){\makebox(0,0){$0$}}%
\put(15.0000,-7.0000){\makebox(0,0){$\bullet$}}%
\put(15.0000,-9.0000){\makebox(0,0){$\bullet$}}%
\put(15.0000,-11.0000){\makebox(0,0){$\vdots$}}%
\put(15.0000,-14.0000){\makebox(0,0){$\bullet$}}%
\put(15.0000,-16.0000){\makebox(0,0){$\bullet$}}%
\put(15.0000,-18.0000){\makebox(0,0){$\bullet$}}%
\put(15.0000,-20.0000){\makebox(0,0){$\bullet$}}%
\put(15.0000,-22.0000){\makebox(0,0){$\bullet$}}%
\put(15.0000,-24.0000){\makebox(0,0){$\bullet$}}%
\put(15.0000,-26.0000){\makebox(0,0){$\bullet$}}%
%
\special{pn 8}%
\special{pa 1400 2800}%
\special{pa 1400 400}%
\special{pa 1400 400}%
\special{pa 1400 400}%
\special{fp}%
\put(18.0000,-7.0000){\makebox(0,0){$\cdots$}}%
\put(18.0000,-9.0000){\makebox(0,0){$\cdots$}}%
\put(18.0000,-11.5000){\makebox(0,0){$\cdots$}}%
\put(18.0000,-14.0000){\makebox(0,0){$\cdots$}}%
\put(18.0000,-16.0000){\makebox(0,0){$\cdots$}}%
\put(18.0000,-18.0000){\makebox(0,0){$\cdots$}}%
\put(18.0000,-20.0000){\makebox(0,0){$\cdots$}}%
\put(18.0000,-22.0000){\makebox(0,0){$\cdots$}}%
\put(18.0000,-24.0000){\makebox(0,0){$\cdots$}}%
\put(18.0000,-26.0000){\makebox(0,0){$\cdots$}}%
\put(21.0000,-26.0000){\makebox(0,0){$*$}}%
\put(21.0000,-24.0000){\makebox(0,0){$\bullet$}}%
\put(21.0000,-22.0000){\makebox(0,0){$\bullet$}}%
\put(21.0000,-20.0000){\makebox(0,0){$\bullet$}}%
\put(21.0000,-18.0000){\makebox(0,0){$\bullet$}}%
\put(21.0000,-16.0000){\makebox(0,0){$\bullet$}}%
\put(21.0000,-14.0000){\makebox(0,0){$\bullet$}}%
\put(21.0000,-11.0000){\makebox(0,0){$\vdots$}}%
\put(21.0000,-9.0000){\makebox(0,0){$\bullet$}}%
\put(21.0000,-7.0000){\makebox(0,0){$\bullet$}}%
\put(23.5000,-24.0000){\makebox(0,0){$\circ$}}%
\put(23.5000,-22.0000){\makebox(0,0){$\bullet$}}%
\put(23.5000,-20.0000){\makebox(0,0){$\bullet$}}%
\put(23.5000,-18.0000){\makebox(0,0){$\bullet$}}%
\put(23.5000,-16.0000){\makebox(0,0){$\bullet$}}%
\put(23.5000,-14.0000){\makebox(0,0){$\bullet$}}%
\put(23.5000,-11.0000){\makebox(0,0){$\vdots$}}%
\put(23.5000,-9.0000){\makebox(0,0){$\bullet$}}%
\put(23.5000,-7.0000){\makebox(0,0){$\bullet$}}%
\put(26.0000,-22.0000){\makebox(0,0){$\bullet$}}%
\put(26.0000,-20.0000){\makebox(0,0){$\bullet$}}%
\put(26.0000,-18.0000){\makebox(0,0){$\bullet$}}%
\put(26.0000,-16.0000){\makebox(0,0){$\bullet$}}%
\put(26.0000,-14.0000){\makebox(0,0){$\bullet$}}%
\put(26.0000,-11.0000){\makebox(0,0){$\vdots$}}%
\put(26.0000,-9.0000){\makebox(0,0){$\bullet$}}%
\put(26.0000,-7.0000){\makebox(0,0){$\bullet$}}%
\put(29.0000,-7.0000){\makebox(0,0){$\cdots$}}%
\put(29.0000,-9.0000){\makebox(0,0){$\cdots$}}%
\put(29.0000,-11.5000){\makebox(0,0){$\cdots$}}%
\put(29.0000,-14.0000){\makebox(0,0){$\cdots$}}%
\put(29.0000,-16.0000){\makebox(0,0){$\cdots$}}%
\put(29.0000,-18.0000){\makebox(0,0){$\cdots$}}%
\put(29.0000,-20.0000){\makebox(0,0){$\cdots$}}%
\put(29.0000,-22.0000){\makebox(0,0){$\cdots$}}%
\put(32.0000,-22.0000){\makebox(0,0){$*$}}%
\put(32.0000,-20.0000){\makebox(0,0){$\bullet$}}%
\put(32.0000,-18.0000){\makebox(0,0){$\bullet$}}%
\put(32.0000,-16.0000){\makebox(0,0){$\bullet$}}%
\put(32.0000,-14.0000){\makebox(0,0){$\bullet$}}%
\put(32.0000,-11.0000){\makebox(0,0){$\vdots$}}%
\put(32.0000,-9.0000){\makebox(0,0){$\bullet$}}%
\put(32.0000,-7.0000){\makebox(0,0){$\bullet$}}%
\put(34.5000,-7.0000){\makebox(0,0){$\bullet$}}%
\put(34.5000,-9.0000){\makebox(0,0){$\bullet$}}%
\put(34.5000,-11.0000){\makebox(0,0){$\vdots$}}%
\put(34.5000,-14.0000){\makebox(0,0){$\bullet$}}%
\put(34.5000,-16.0000){\makebox(0,0){$\bullet$}}%
\put(34.5000,-18.0000){\makebox(0,0){$\circ$}}%
\put(34.5000,-20.0000){\makebox(0,0){$\circ$}}%
\put(37.0000,-18.0000){\makebox(0,0){$\circ$}}%
\put(37.0000,-16.0000){\makebox(0,0){$\bullet$}}%
\put(37.0000,-14.0000){\makebox(0,0){$\bullet$}}%
\put(37.0000,-11.0000){\makebox(0,0){$\vdots$}}%
\put(37.0000,-9.0000){\makebox(0,0){$\bullet$}}%
\put(37.0000,-7.0000){\makebox(0,0){$\bullet$}}%
\put(39.5000,-7.0000){\makebox(0,0){$\bullet$}}%
\put(39.5000,-9.0000){\makebox(0,0){$\bullet$}}%
\put(39.5000,-11.0000){\makebox(0,0){$\vdots$}}%
\put(39.5000,-14.0000){\makebox(0,0){$\bullet$}}%
\put(39.5000,-16.0000){\makebox(0,0){$\bullet$}}%
\put(42.0000,-16.0000){\makebox(0,0){$*$}}%
\put(42.0000,-14.0000){\makebox(0,0){$\bullet$}}%
\put(42.0000,-11.0000){\makebox(0,0){$\vdots$}}%
\put(42.0000,-9.0000){\makebox(0,0){$\bullet$}}%
\put(42.0000,-7.0000){\makebox(0,0){$\bullet$}}%
\put(44.5000,-7.0000){\makebox(0,0){$\bullet$}}%
\put(44.5000,-9.0000){\makebox(0,0){$\circ$}}%
\put(47.0000,-7.0000){\makebox(0,0){$*$}}%
\put(47.0000,-5.0000){\makebox(0,0){$i_4$}}%
\put(42.0000,-5.0000){\makebox(0,0){$i_3$}}%
\put(32.0000,-5.0000){\makebox(0,0){$i_2$}}%
\put(21.2000,-5.0000){\makebox(0,0){$i_1$}}%
%
\special{pn 8}%
\special{pa 4800 600}%
\special{pa 4800 800}%
\special{pa 4350 800}%
\special{pa 4350 1700}%
\special{pa 3350 1700}%
\special{pa 3350 2300}%
\special{pa 2250 2300}%
\special{pa 2250 2700}%
\special{pa 1400 2700}%
\special{pa 1400 2700}%
\special{pa 1400 2700}%
\special{fp}%
%
\special{pn 8}%
\special{pa 4900 600}%
\special{pa 1100 600}%
\special{pa 1100 600}%
\special{pa 1100 600}%
\special{fp}%
\put(13.0000,-7.0000){\makebox(0,0){$j_4$}}%
\put(13.0000,-16.0000){\makebox(0,0){$j_3$}}%
\put(13.0000,-22.0000){\makebox(0,0){$j_2$}}%
\put(13.0000,-26.0000){\makebox(0,0){$j_1$}}%
\end{picture}%
\end{center}
($*$ denotes an extremal point of $A$ and $\circ$ denotes a point $(i,j) \in \mathrm{conv}(A) \setminus A$.)
\medskip

Let $I \subset S$ be a graded ideal and $H$ the Hilbert function 
of $S/I$.
Then there exists the smallest extremal area $A$
such that $I$ admits $A$.
Although, there is not always a graded ideal $L \in \mathcal{L}_H^A$
which has the maximal graded Betti numbers among graded ideals in $\mathcal{L}_H^A$,
Theorem \ref{laststand} guarantees the existence of such an ideal
if we extend $A$ to $\mathrm{conv}(A)$.
Also,
since $A \subset Q$ implies $\mathrm{conv}(A) \subset Q$,
we can consider the same property for squarefree monomial
ideals by using the squarefree operation.
\bigskip

\noindent
\textbf{Acknowledgement:}
I would like to thank the referee for pointing out the relation between the results in \S 4
and \cite[Corollary 3.2]{Bj}.

\end{document}